\documentclass[11pt]{amsart}
\usepackage{amscd}      
\usepackage{amsbsy, amsthm, eucal, graphics,epic,amsbsy, stmaryrd, epic, eepic}
\usepackage{amssymb}
\usepackage{amsmath}
\usepackage[all]{xy}

\usepackage{latexsym}
\usepackage{epsfig}
\usepackage{amsfonts}
\usepackage{enumerate}

\newtheorem{prop}{Proposition}[section]
\newtheorem{lem}[prop]{Lemma}

\newtheorem{cor}[prop]{Corollary}
\newtheorem{thm}[prop]{Theorem}
\newtheorem{rmk}[prop]{Remark}

\newtheorem{defn}[prop]{Definition}

\newenvironment{pf}{\begin{trivlist}\item[]{\sc Proof.}}%
            {\nolinebreak $\Box$ \end{trivlist}}

\newcommand{\noprint}[1]{}

\renewcommand{\tilde}{\widetilde}

\newcommand{\YY}{{\mathfrak Y}}

\newcommand{\MM}{{\mathfrak M}}

\newcommand{\CC}{{\mathfrak C}}

\newcommand{\cc}{{\mathbb C}}

\newcommand{\ldiag}[1]%
       {\makebox[0cm]{${\scriptstyle#1}\downarrow\phantom{\scriptstyle#1}$}}
\newcommand{\ldiagup}[1]%
       {\makebox[0cm]{${\scriptstyle#1}\uparrow\phantom{\scriptstyle#1}$}}
\newcommand{\rdiag}[1]%
       {\makebox[0cm]{$\phantom{\scriptstyle#1}\downarrow{\scriptstyle#1}$}}
\newcommand{\sediagr}[1]%
       {\makebox[0cm]{$\phantom{\scriptstyle#1}\searrow{\scriptstyle#1}$}}
\newcommand{\nediagr}[1]%
       {\makebox[0cm]{$\phantom{\scriptstyle#1}\nearrow{\scriptstyle#1}$}}
\newcommand{\rdiagup}[1]%
       {\makebox[0cm]{$\phantom{\scriptstyle#1}\uparrow{\scriptstyle#1}$}}
\newcommand{\swdiag}[1]%
       {\makebox[0cm]{$\phantom{\scriptstyle#1}\swarrow{\scriptstyle#1}$}}
\newcommand{\sediag}[1]%
       {\makebox[0cm]{${\scriptstyle#1}\searrow\phantom{\scriptstyle#1}$}}
\newcommand{\nediag}[1]%
       {\makebox[0cm]{${\scriptstyle#1}\nearrow\phantom{\scriptstyle#1}$}}

\newcommand{\doublearrowstack}[2]%
                      {{{{\scriptstyle#1}\atop{\textstyle\longrightarrow}}\atop{{\textstyle\longrightarrow}\atop{\scriptstyle#2}}}}
\newcommand{\rightleftarrowstack}[2]%
                      {{{{\scriptstyle#1}\atop{\textstyle\longrightarrow}}\atop{{\textstyle\longleftarrow}\atop{\scriptstyle#2}}}}
\newcommand{\leftrightarrowstack}[2]%
                      {{{{\scriptstyle#1}\atop{\textstyle\longleftarrow}}\atop{{\textstyle\longrightarrow}\atop{\scriptstyle#2}}}}

\newcommand{\overtoparrow}%
{\makebox[0cm]{\beginpicture \setcoordinatesystem units
<.8cm,.4cm> point at 0 0 \setplotarea x from -3 to 3, y from 0 to
1 \setquadratic \plot -3 0 0 1 3 0 / \put{\vector(3,-1){0}}[Bl] at
3 0
\endpicture}}

\newcommand{\underbottomarrow}%
{\makebox[0cm]{\beginpicture \setcoordinatesystem units
<.8cm,.4cm> point at 0 0 \setplotarea x from -3 to 3, y from 0 to
1 \setquadratic \plot -3 1 0 0 3 1 / \put{\vector(3,1){0}}[Bl] at
3 1
\endpicture}}

\newcommand{\ses}[5]%
{0\longrightarrow#1\stackrel{#2}{ \longrightarrow}#3\stackrel{#4}{
\longrightarrow}#5\longrightarrow0}

\newcommand{\dt}[6]%
{#1\stackrel{#2}{longrightarrow}#3
\stackrel{#4}{\longrightarrow}#5 \stackrel{#6}{\longrightarrow}
#1[1]}

\newcommand{\cat}[1]%
{(\mbox{\rm #1})}

\def\bbA{{\mathbb A}}  \def\bbC{{\mathbb C}}  

\def\bbG{{\mathbb G}}    
\def\bbQ{{\mathbb Q}} \def\bbZ{{\mathbb Z}}   \def\bbN{{\mathbb N}}

\newcommand{\clD}{{\mathcal{D}}}
\newcommand{\clC}{{\mathcal{C}}}

\newcommand{\clG}{{\mathcal{G}}}
\newcommand{\clF}{{\mathcal F}}
\newcommand{\clK}{{\mathcal{K}}}
\newcommand{\clL}{{\mathcal{L}}}
\newcommand{\clM}{{\mathcal M}}
\newcommand{\clN}{{\mathcal N}}

\newcommand{\clP}{{\mathcal P}}
\newcommand{\clT}{{\mathcal{T}}}

\newcommand{\clO}{{\mathcal{O}}}

\newcommand{\clX}{{\mathcal{X}}}
\def\longto{\longrightarrow}
\def\isomto{\stackrel{\sim}{\longto}}
\providecommand{\abs}[1]{\lvert#1\rvert}
\newcommand{\Spec}{\operatorname{Spec}}
\newcommand{\Pic}{\operatorname{Pic}}

\title[Root gerbes I: structure of genus $0$ moduli spaces]{Gromov-Witten theory of root gerbes I:\\ structure of genus $0$ moduli spaces}
\date{\today}
\author{Elena Andreini}
\address{SISSA\\via Bonomea, 265\\34136 Trieste TS\\Italy}
\email{andreini.elena@gmail.com}
\author{Yunfeng Jiang}
\address{Department of Mathematics\\ Imperial College London\\ South Kensington Campus\\ London SW7 2AZ\\ United Kingdom}
\email{y.jiang@imperial.ac.uk}
\author{Hsian-Hua Tseng}
\address{Department of Mathematics\\ Ohio State University\\ 100 Math Tower, 231 West 18th Ave.\\Columbus\\ OH 43210\\ USA}
\email{hhtseng@math.ohio-state.edu}

\begin{document}
\begin{abstract}
Let $X$ be a smooth complex projective algebraic variety. Given a line bundle $\mathcal{L}$ over $X$ and an integer $r>1$ one defines the stack $\sqrt[r]{\mathcal{L}/X}$ of $r$-th roots of $\mathcal{L}$. Motivated by Gromov-Witten theoretic questions, in this paper we analyze the structure of moduli stacks of genus $0$ twisted stable maps to $\sqrt[r]{\mathcal{L}/X}$. Our main results are explicit constructions of moduli stacks of genus $0$ twisted stable maps to $\sqrt[r]{\mathcal{L}/X}$ starting from moduli stack of genus $0$ stable maps to $X$. As a consequence, we prove an exact formula expressing genus $0$ Gromov-Witten invariants of $\sqrt[r]{\mathcal{L}/X}$ in terms of those of $X$. 
\end{abstract}

\bibliographystyle{alpha}
\sloppy \maketitle
\tableofcontents

\section{Introduction}
Orbifold Gromov-Witten theory, constructed in symplectic category by Chen-Ruan \cite{CR2}  and in algebraic category by Abramovich, Graber and Vistoli \cite{AGV1}, \cite{AGV2}, has been an area of active research in recent years. Calculations of orbifold Gromov-Witten invariants in examples present numerous new challenges, see \cite{CCLT}, \cite{CCIT}, \cite{ts}, and \cite{BC} for examples. 

\'Etale gerbes over a smooth base provide interesting examples of smooth Deligne-Mumford stacks. Let $\clX$ be a smooth Deligne-Mumford stack and  $G$  a finite group  scheme over $\clX$. Intuitively one can think of a $G$-banded gerbe over  $\clX$ as a fibre bundle over $\clX$  with fibre the classifying stack $BG$. A detailed definition of gerbes can be found in, for example, \cite{Gir}, \cite{Bree94}, \cite{ehkv}. We are interested in computing Gromov-Witten theory of $G$-banded gerbes. 

Physics considerations have suggested that the geometry of \'etale gerbes possesses certain very intriguing structure. The so-called {\em decomposition conjecture} \cite{HHPS} in physics maybe interpreted mathematically as a philosophy saying that the geometry of an \'etale gerbe is equivalent to the geometry of certain disconnected space twisted by a $U(1)$-gerbe. In-depth discussions on various mathematical aspects of this conjecture can be found in \cite{TT}. 

The Gromov-Witten theoretic version of the decomposition conjecture, which can be formulated for arbitrary $G$-gerbes more general than $G$-banded gerbes, states that Gromov-Witten theory of the $G$-gerbe is equivalent to certain twist of the Gromov-Witten theory of some \'etale cover of the base. A detailed discussion of the conjecture in full generality can be found in \cite{TT}. For $G$-banded gerbes this conjecture states that the Gromov-Witten theory of a $G$-banded gerbe over $\clX$ is equivalent to (certain twists of) the Gromov-Witten theory of the disjoint union of $|\text{Conj}(G)|$ copies\footnote{Here $\text{Conj}(G)$ is the set of conjugacy classes of $G$.} of $\clX$ after a change of variables. Computations of Gromov-Witten invariants of \'etale gerbes is thus intimated connected to the decomposition conjecture.

The simplest examples of $G$-gerbes are {\em trivial gerbes}. The trivial $G$-gerbe over a Deligne-Mumford stack $\clX$ is the product $\clX\times BG$. In \cite{AJT} the computation of Gromov-Witten invariants of $\clX\times BG$ is handled as a special case of a general product formula for orbifold Gromov-Witten invariants of product Deligne-Mumford stacks $\clX\times\YY$. As a consequence the decomposition conjecture is proven for trivial $G$-gerbes. 

An interesting class of non-trivial gerbes is provided by {\em root gerbes} associated to line bundles. This is the first of two papers in which we study Gromov-Witten theory of root gerbes of line bundles over smooth projective varieties, with the decomposition conjecture in mind. The present paper is devoted to study the genus $0$ Gromov-Witten theory of root gerbes. 

Let $X$ be a smooth complex projective variety and $\mathcal{L}\to X$ a line bundle. Given an integer $r>0$, let $$\clG:=\sqrt[r]{\mathcal{L}/X}\to X$$ be the stack of $r$-th roots of $\mathcal{L}$ over $X$. It can be shown that $\clG\to X$ is a $\mu_r$-banded gerbe over $X$. Such a gerbe is called a root gerbe. Constructions and properties of root gerbes are briefly reviewed in Section \ref{root-gerbe}. In order to study the Gromov-Witten theory we consider moduli spaces $\mathcal{K}_{0,n}(\clG,\beta)$ of genus $0$ twisted stable maps to $\clG$. By composing a twisted stable map to $\clG$ with the structure map $\clG\to X$, one can define a morphism 
\begin{equation}\label{natural_map}
\mathcal{K}_{0,n}(\clG,\beta)\to \overline{M}_{0,n}(X,\beta),
\end{equation}
 where $\overline{M}_{0,n}(X,\beta)$ is a moduli space of  genus $0$ stable maps to $X$. The main idea used in our approach to Gromov-Witten theory of root gerbes is to compare Gromov-Witten invariants of $\clG$ with Gromov-Witten invariants of the base $X$ using the morphism (\ref{natural_map}). In the present paper, this idea is realized by our main results, Theorems \ref{mod_space_to_gerbe_is_gerbe} and \ref{mod_space_to_root_gerbe_is_root_gerbe}, on the structures of the moduli spaces $\mathcal{K}_{0,n}(\clG,\beta)$. Roughly speaking, these structure results state that components of $\clK_{0,n}(\clG, \beta)$ are $\mu_r$-gerbes over certain base stacks constructed from $\overline{M}_{0,n}(X, \beta)$ using log geometry. More details can be found in Sections \ref{section:gerbe_structure} and \ref{section:root_gerbe_structure}. Our results extends a result of \cite{BC} for the gerbe $B\mu_r$.  Our proofs are based on a detailed analysis of the moduli spaces $\mathcal{K}_{0,n}(\clG, \beta)$, and use heavily the results of \cite{MO} and \cite{OLogCurv}.

As a consequence of our main structure results Theorems \ref{mod_space_to_gerbe_is_gerbe} and \ref{mod_space_to_root_gerbe_is_root_gerbe}, we prove a comparison result between virtual fundamental classes of $\mathcal{K}_{0,n}(\clG, \beta)$ and $\overline{M}_{0,n}(X,\beta)$, see Theorem \ref{pushforward_vir}. This comparison result yields an explicit computation of genus $0$ Gromov-Witten invariants of $\clG$ in terms of genus $0$ Gromov-Witten invariants of $X$, which is Theorem \ref{GW-inv1}. A reformulation of Theorem \ref{GW-inv1} in terms of generating functions confirms the decomposition conjecture for genus $0$ Gromov-Witten theory of $\clG$, see Theorem \ref{decomp_genus_0}.

The paper is organized as follows. Section \ref{pre} contains discussions on some preparatory materials. In Section \ref{moduli_maps_root_gerbe} we carry out the needed analysis on the structure of the moduli spaces of twisted stable maps to root gerbes. In Section \ref{GW_theory} we prove results on virtual fundamental classes and Gromov-Witten invariants, in particular the decomposition conjecture in genus $0$. In Appendix \ref{banded_abelian_gerbes_case} we discuss extensions of our results to banded abelian gerbes.

\subsection*{Conventions}
Unless otherwise mentioned, we work over $\bbC$ throughout this paper. By an {\em algebraic stack} we mean an  algebraic stack over $\bbC$ in the sense of \cite{Art74}. By a {\em Deligne-Mumford stack} we mean an algebraic stack over $\bbC$ in the sense of \cite{DM69}. We assume moreover all stacks (and schemes) are  quasi-separated, locally noetherian, locally of finite type. From time to time we use the notation $x\in X$ to indicate that $x$ is a geometric point of $X$. Following \cite{KatoLog}, logarithmic structures are considered on the \'etale site of schemes. For the  extension of logarithmic structures to stacks,  see \cite{OlssLog03}. Given a scheme (or a stack) $X$, a geometric point $x$ of $X$, and a sheaf of sets $\clF$ on $X$, according to the standard notation  we denote by $\clF_{\overline{x}}$ the stalk of $\clF$ at $x$ in the \'etale topology.  A gerbe  is an algebraic stack as in \cite{l-mb} Definition 3.15.

\subsection*{Acknowledgments} We thank D. Abramovich, A. Bayer,  K. Behrend,  B. Fantechi,  P. Johnson, A. Kresch, F. Nironi,  E. Sharpe,  Y. Ruan and A. Vistoli for valuable discussions. H.-H. T. is grateful to T. Coates, A. Corti, H. Iritani, and X. Tang for related collaborations. Y. J. and H.-H. T. thanks Mathematical Sciences Research Institute for hospitality and support of a visit in spring 2009 during which part of this paper was written. H.-H. T. is supported in part by NSF grant DMS-0757722.

\section{Preliminaries}\label{pre}

\subsection{Twisted stable maps}
We recall the definition of twisted curve here, see \cite{AGV1}, \cite{AGV2}, \cite{AV} for more details.
\begin{defn}[\cite{AV}, Definition 4.1.2]
A twisted nodal $n$-pointed curve over a scheme $S$ is a morphism $\mathcal{C}\to S$ together with $n$ closed substacks $\sigma_i\subset \mathcal{C}$ such that
\begin{itemize}
\item
$\mathcal{C}$ is a tame Deligne-Mumford stack, proper over $S$, and \'etale locally is a nodal curve over $S$;

\item
$\sigma_i\subset \mathcal{C}$ are disjoint closed substacks in the smooth locus of $\mathcal{C}\to S$;

\item
$\sigma_i\to S$ are \'etale gerbes;

\item
the map $\mathcal{C}\to C$ to the coarse moduli space $C$ is an isomorphism away from marked points and nodes.
\end{itemize}
\end{defn}

By definition the genus of a twisted curve $\mathcal{C}\to S$ is the genus of its coarse moduli space $C\to S$. 

Throughout this paper we will always assume that twisted curves are {\em balanced}, i.e. at any twisted node, the local group acts on the two branches by opposite characters.

Let $S$ be a noetherian  scheme and let  $\mathcal{X}/S$ be a proper  Deligne-Mumford stack over $S$ with projective coarse moduli space $X\to S$.
We fix an ample invertible sheaf $\clO_X(1)$ over $X$. 
 Let $\mathcal{K}_{g,n}(\mathcal{X}, \beta)$  be the fibered
category over $S$ which to any $S$-scheme $T$ associates the groupoid of the following data:
\begin{itemize}
\item
A twisted $n$-pointed curve $(\mathcal{C}/T, \{\sigma_{i}\})$ over $T$;
\item
A representable morphism $f : \mathcal{C}\to \mathcal{X}$ such that the induced morphism $\bar{f}: C\to  X$ between coarse moduli spaces is an $n$-pointed stable map of degree $\beta\in H^+_2(X, \mathbb{Z})$ (i.e. $\bar{f}_*[C]=\beta$).
\end{itemize}
According to \cite{AV}, Theorem 1.4.1 the fibered category $\mathcal{K}_{g,n}(\mathcal{X}, \beta)$ is a Deligne-Mumford stack  proper  over $S$. 

As discussed in \cite{AGV2}, there exist  evaluation maps: 
$$ev_{i}: \mathcal{K}_{0,n}(\mathcal{X},\beta)\to \bar{I}(\mathcal{X}), \quad 1\leq i\leq n$$ taking values in the {\em rigidified inertia stack} $\bar{I}(\mathcal{X})$ of $\mathcal{X}$. This map is obtained as follows. The rigidified inertia stack $\bar{I}(\clX)$ may be defined as the stack of cyclotomic gerbes in $\clX$, i.e. representable morphisms from cyclotomic gerbes to $\clX$. The evaluation map $ev_i$ is defined to map a twisted stable map $f: (\mathcal{C}/T, \{\sigma_{i}\})\to \mathcal{X}$ to its restriction to the $i$-th marked gerbe,  
$$f|_{\sigma_i}: \sigma_i\to \mathcal{X},$$
which is an object of $\bar{I}(\mathcal{X})$. 

The rigidified inertia stack $\bar{I}(\clX)$ has an alternative description. Define the {\em inertia stack} of $\clX$ to be the fiber product over the diagonal: $$I\clX:= \clX\times_{\clX\times_S \clX}\clX.$$ 
By definition, objects of $I\clX$ are pairs $(x, g)$ where $x$ is an object of $\clX$ and $g$ is an element of the automorphism group of $x$. The rigidified inertia stack $\bar{I}(\clX)$ is obtained from $I\clX$ by applying the rigidification procedure (\cite{ACV03}, \cite{AOV07}). More details can be found in e.g. \cite{AGV1}.

\subsection{Root gerbes}\label{root-gerbe}

We recall the notion of root gerbes. Let $X$ be a smooth projective variety and let $\clL$ be  a line bundle over $X$ corresponding to a morphism $\phi_\clL: X\to B\bbC^*$.
For an integer $r>0$ let $\theta_r: B\bbC^*\to B\bbC^*$ be the morphism  induced  by  the  $r$-th power homomorphism  $\bbC^*\stackrel{(\cdot)^r}{\to} \bbC^*$. The composite morphism $\theta_r\circ\phi_\clL: X \to B\bbC^*$ corresponds to $\clL^{\otimes r}$.

\begin{defn}
The stack $\sqrt[r]{\clL/X}$ of $r$-th roots of $\clL$ is  defined as
$$
\sqrt[r]{\clL/X}:= X\times_{\phi_\clL, B\bbC^*,\theta_r} \mbox{B}\bbC^*.
$$
Explicitly it can be described as the $X$-groupoid whose objects
over $(Y,f:Y\to X)$ are pairs $(M,\varphi)$, with $M$ a line bundle over $Y$ and $\phi: M^{\otimes r}\to f^*\clL$ an
isomorphism. An arrow from $(M,\varphi)$ to $(N,\psi)$ lying over a $X$-morphism $h: (Y,f)\to (Z,g)$ is an isomorphism
$\rho: M\to h^{*}N$ such that
$\varphi$ fits in the following commutative diagram
\[\xymatrix{
M^{\otimes r} \ar[r]^{\rho^{\otimes r}} \ar[d]_{\varphi} & h^{*}N^{\otimes r} \ar[d]^{h^{*}\psi} \\
f^{*}\clL  \ar[r]_{\cong}  & h^{*}g^{*}\clL, }
\]
where the bottom arrow is the canonical isomorphism.
\end{defn}

The following proposition follows easily from the definition.

\begin{prop}\label{root_stack_quotient}
The stack $\sqrt[r]{\clL/X}$ is the quotient stack $[\clL^{\times}/\mathbb{C}^*]$, where $\clL^{\times}$ is the principal $\mathbb{C}^*$-bundle obtained by deleting the zero section of $\clL$, and $\mathbb{C}^*$ acts on $\clL^\times$ via $\lambda \cdot z= \lambda^r z, \lambda \in \mathbb{C}^*, z\in \clL^\times$. In particular $\sqrt[r]{\clL/X}$ is a Deligne-Mumford stack.
\end{prop}

\begin{pf}
It is enough to observe that the following  diagram is 2-cartesian
\begin{eqnarray*}
\xymatrix{
\clL^{\times}\ar[d]\ar[r]\ar@{}[rd]|{\square} & \sqrt[r]{\clL/X}\ar[r]\ar[d]\ar@{}[rd]|{\square}& X\ar[d]\\
pt\ar[r]& B\bbC^*\ar[r]^{\theta_r} & B\bbC^*.  
}
\end{eqnarray*}
\end{pf}


\begin{rmk}
The morphism $\theta_r:\mbox{B}\bbC^* \to \mbox{B}\bbC^*$ is a $\mu_r$-gerbe, because of the Kummer exact sequence
$$1\to \mu_{r}\to \mathbb{C}^*\stackrel{(\cdot)^{r}}{\to} \mathbb{C}^*\to 1.$$
Hence   $\sqrt[r]{\clL/X}\to X$ is a $\mu_r$-gerbe.
\end{rmk}

\begin{rmk}
The stack $\sqrt[r]{\clL/X}$ may also be constructed as a toric stack bundle \cite{Jiang}. 
\end{rmk}
It is also  possible to take roots of ``line bundles with sections''. Let $\clL$ be a line bundle over $X$ and let $\sigma$ be a section of $\clL$. The data $(\clL,\sigma)$ correspond to a morphism $\phi_{\clL,\sigma}: X\to [\bbA^1/\bbC^*]$. Let $\theta_r:  [\bbA^1/\bbC^*]\to  [\bbA^1/\bbC^*]$ be the morphism induced by the $r$-th power morphisms on $\bbA^1$ and  $\bbC^*$. The morphism $\theta_r\circ \phi_{\clL,\sigma}$ corresponds to the pair  $(\clL^{\otimes r}, \sigma^r)$.
 The stack $\sqrt[r]{(\clL,\sigma)/X}$ of $r$-th roots of $\clL$ with the section $\sigma$ is defined as 
$$ 
\sqrt[r]{(\clL,\sigma)/X}:= [\mathbb{A}^1/\bbC^*]\times_{ [\mathbb{A}^1/\bbC^*],\theta_r} [\mathbb{A}^1/\bbC^*].
$$
The stack constructed in this way is isomorphic to $X$ outside the vanishing locus $Z(\sigma)\subset X$ of $\sigma$, while the reduced substack of the closed substack mapping to  $Z(\sigma)$ is a $\mu_r$-gerbe over $Z(\sigma)$. 
Note that given a divisor $D\subset X$  there is an associated line bundle with 
a canonical section which vanishs on $D$. Therefore in the following we will also talk about roots of divisors. 

\subsection{Line bundles over twisted curves}\label{Pic-group-tw-curve-subs}
We recall some results about line bundles over twisted curves. In \cite{Ca} there is an explicit description of the Picard group of a smooth  twisted curve. Let $\clC$ be a smooth twisted curve over $\text{Spec}\, \bbC$. Let $C$ be the coarse curve and $D_i\in C, 1\leq i\leq n$ the marked points. It is known that $\clC$ can be constructed from its coarse curve $C$ by applying the $r_i$-th root construction to the divisor $D_i$, for all $1\leq i\leq n$. (Here $r_i\in \bbN$.) Let $\clT_i, 1\leq i\leq n$ be the tautological line bundles associated by the root construction and $\tau_i, 1\leq i\leq n$ their tautological sections.
\begin{lem}[\cite{Ca}, Corollary 2.12]
\label{coroPic_clC}
Let $\clL$ be an invertible sheaf on $\clC$. Then there exists an invertible sheaf $L$ on $C$ and integers $k_i$ satisfying $0\leq k_i\leq r_i-1$ such that
\begin{equation*}
\clL\simeq \pi^* L\otimes \prod_{i=1}^{n}\clT^{k_i}_i.
\end{equation*}
Moreover the integers $k_i$ are unique, $L$ is unique up to isomorphism. 
\end{lem}
There is an analogous description for the global sections of invertible sheaves on $\clC$. 
\begin{lem}[\cite{Ca}, Corollary 2.13]
Given the decomposition in Lemma \ref{coroPic_clC}, every global section of $\clL$ is of the form $\pi^* s\otimes \tau_1^{k_1}...\otimes \tau_n^{k_n}$ for a unique global section $s$ of $L$, where $\tau_i$ is the tautological section of $\clT_i$.
\end{lem}

Lemma \ref{coroPic_clC} can be rephrased as saying that $\Pic{\clC}$ is an extension of $\Pic{C}$ by a finite abelian group, namely
\begin{equation*}
 1\to \Pic{C}\to \Pic{\clC}\to \oplus_{i=1}^n \bbZ_{r_i}\to 1,
\end{equation*}
 where $r_i$ are the orders of the stabilizers of stack points. 
\begin{rmk}
The same description of $\Pic{\clC}$ holds when $\clC$ is not smooth but  has only untwisted nodes. 
\end{rmk}

The Picard groups of nodal  twisted curves over $\Spec{\bbC}$ admit a similar description.
This is shown in e.g. \cite{Ch06}. We sketch the argument for the reader's convenience.
\begin{lem}[See \cite{Ch06}, Theorem 3.2.3]\label{Pic_nod_tw_lemm}
Let $\clC$ be an unmarked twisted curve with nodes $e_1,..,e_s$. Let   $\gamma_j$ be the  order of the stabilizer of the node $e_j$.
Then the following exact sequence holds:
\begin{equation*}
1\to\Pic{C}\to\Pic{\clC}\to \prod_{j=1}^s\bbZ/\gamma_j\bbZ\to 1.
\end{equation*}
\end{lem}
\begin{pf}
Let $\pi: \clC\to C$ be the map to the coarse curve. Consider the exact sequence of complexes over $C$ given by
\begin{equation*}
1\to\pi_*\bbG_m\to R\pi_*\bbG_m\to R\pi_*\bbG_m/\pi_*\bbG_m\to 1.  
\end{equation*}
Notice that $\pi_*\mu_r=\mu_r$ and $\pi_*\bbG_m=\bbG_m$. Therefore they are complexes concentrated in
degree zero. The long Hypercohomology  exact sequence gives
\begin{equation*}
1\to H^1(C,\bbG_m)\to H^1(\clC,\bbG_m)\to H^1(R\pi_*\bbG_m/\pi_*\bbG_m)\to 1.
\end{equation*} 
This sequence is exact on the left because $E^{p,q}_2:=H^p(C,H^q(R\pi_*\bbG_m/\bbG_m)) 
$ abuts to $\mathbb{H}^{p+q}(C,R\pi_*\bbG_m/\bbG_m)$. The sheaf $H^q(R\pi_*\bbG_m/\bbG_m)$ is equal to $R^q\pi_*\bbG_m$ and does not vanish for $q>0$. By \cite{ACV03}, Proposition A.0.1, the stalk of $R^q\pi_*\bbG_m$ is canonically isomorphic to $H^q(Aut(p),\bbG_{m,p})$ where $p$ is a geometric point of $C$.
This sequence is exact on the right because $H^2(C,\bbG_m)=0$ for $C$ a genus zero nodal curve.
The result follows by observing that $H^q(\mu_r,\bbG_m)=\bbZ/r\bbZ$ for $q$ odd and is trivial for $q$ even.
\end{pf}
\begin{rmk}
The above proof generalizes to nodal marked twisted curves.
\end{rmk}
\subsubsection{Normalization of  twisted curves}

It is very useful to describe a twisted stable map over a point   $\tilde{f}:\clC\to\clG$  in terms of the induced morphism $\tilde{f}\circ\nu:\tilde{\clC}\to \clG$, where $\tilde{\clC}$ is the normalization of $\clC$. This morphism is still a twisted stable map (with possibly disconnected domain). According to  \cite{Vist89} the normalization  of a reduced stack  $\clX$ is defined in the following way. Let $R\rightrightarrows U$ be a presentation of $\clX$. Let $\tilde{R}$ and $\tilde{U}$ be the normalizations of $R$ and $U$. It is possible to lift the structure morphisms of the groupoid $R\rightrightarrows U$  in such a way that  $\tilde{R}\rightrightarrows \tilde{U}$ is also a groupoid. Moreover the diagonal $\tilde{R}\to \tilde{U}\times \tilde{U}$ is separated and quasi compact. Therefore the groupoid defines an algebraic stack, which is the normalization of $\clX$. In particular the normalization morphism $\nu:\tilde{\clX}\to \clX$ is  representable.

Smooth twisted curves admit line bundles whose fibers carry faithful  representations of the stabilizer groups of the points in the special locus. Those are the tautological line bundles obtained from root constructions. Singular twisted curves over a point also admit line bundles with fibers carrying faithful representations of the stabilizer group  of the nodes. This is the content of Lemma \ref{Pic_nod_tw_lemm}. In this case it is easy to describe those line bundles in terms of tautological line bundles on the normalization of the curve. Assume without loss of generality that $\clC$ is a nodal twisted curve with only one node $\mathcal{E}$ of order $\gamma$. Let $e$ be the image of the node in the coarse moduli space $C$. 
 We have the following   commutative  diagram
\begin{equation*}
\xymatrix{
& 1 &   & & \\
1\ar[r] & \Pic{\tilde{C}}\ar[u]\ar[r] & \Pic{\tilde{\clC}}\ar[r] & \oplus_{i=1}^2 \langle \clT_i\rangle\ar[r] & 1\\
1\ar[r] & \Pic{C}\ar[r]\ar[u]^{\nu^*} & \Pic{\clC}\ar[u]^{\tilde{\nu}^*}\ar[r]& \Pic{\mathcal{E}}\ar[r] \ar[u]& 1\\
1\ar[r] & \tilde{\clO}_e^*/\clO_e^*\ar[u]\ar[r] & \tilde{\clO}^*_{\mathcal{E}}/\clO^*_{\mathcal{E}}\ar[r] \ar[u] & 0\ar[u]  \\
& 1\ar[u] & 1\ar[u] & 
}
\end{equation*}
where $\clO_e$, resp. $\clO_{\mathcal{E}}$, is the local ring at the node $e$, resp. at the twisted node $\mathcal{E}$, and  $\tilde{\clO}_e$, resp. $\tilde{\clO}_{\mathcal{E}}$,    is its integral closure.  Note that $\mathcal{E}\simeq B\mu_{\gamma}$.
Here $\langle \clT_i\rangle$  is the group generated by $\clT_i$ under tensor products.
The  line bundle carrying  a representation of the stabilizers group of the node   corresponding 
to an element $\zeta^k$ of $\mu_\gamma$, where $\zeta$ is the standard generator,  is mapped by the   pullback along the normalization morphism $\tilde{\nu}:\tilde{\clC}\to \clC$ to the pair of line bundles $(\clT_+^k,\clT^{-k}_-)$, where
$\clT_+$, $\clT_-$  are  the tautological line bundles associated to the preimages of the node in the normalization. 
\subsection{Logarithmic geometry and twisted curves}\label{subsection:log-geom-and-tw-curves}
We recall here some basic facts about   logarithmic geometry, which is the natural languge to describe twisted curves. We will use logarithmic geometry to construct the auxiliary stack $\YY_{0,n,\beta}^{\vec{g}}$ in Section \ref{section:stack_YY}.\\
Logarithmic structures have been introduced by Fontaine and Illusie and further studied by Kato \cite{KatoLog}.  A generalization to algebraic stacks can be found in \cite{OlssLog03}. We will consider log structures  on the  \'etale site of  schemes and on the Lisse-\'Etale site (\cite{l-mb} 12.1.2 (i)) of algebraic stacks (see  \cite{OlssLog03}, Definition 5.1).

 Given a scheme $X$, a {\em pre-logarithmic structure}, often called {\em pre-log structure},  consists of a  sheaf of monoids $M$ endowed with a morphism of monoids $\alpha:M\to \clO_X$, where    the structure sheaf is considered as a monoid with the multiplicative structure. Given a monoid or a sheaf of monoids $M$, we denote by $M^*$ the submonoid or the subsheaf of invertible elements.

 When the natural morphism $\alpha^{-1}(\clO_X^*)\to M^*$ is an isomorphism, a pre-log structure is called a {\em log structure}. The quotient $M/\alpha^{-1}(\clO^*_X)$ is usually denoted by $\overline{M}$, and  called the {\em characteristic} or the {\em ghost sheaf}.   There is a canonical way to associate a log structure to a pre-log structure. Given a pre-log structure  $\alpha:M\to \clO_X$, the associated log structure, denoted $M^a$, is defined as the pushout in the category of sheaves of monoids as in the following diagram
\begin{eqnarray}
\xymatrix{
\alpha^{-1}(\clO^*_X)\ar[r]\ar[d]_{\alpha} & M\ar[d]\\
\clO^*_X\ar[r] & M^a.
}\nonumber
\end{eqnarray}
The morphism to the sucture sheaf $\alpha^a:M^a\to\clO_X$ is induced by the pair of morphisms $(\alpha,\iota)$, where  $\iota:\clO^*_X\hookrightarrow\clO_X$ is the canonical inclusion. A scheme endowed with  a log structure $(X,M_X)$ is called a {\em log scheme}. Log schemes form a category.  A morphism between two log schemes $(f,f^b): (X,M_X)\to (Y,M_Y)$
 is a pair consisting of a morphism of schemes $f:X\to Y$ and a morphism of sheaves of monoids  $f^{b}:f^*M_Y\to M_X$  compatible with the morphisms to the  structure sheaf. The pullback of a log structure is defined as the log structure associated to the pre-log structure obtained by taking the inverse image.

A log structure $\clM_X$ over $X$ is called {\em locally free} if for any geometric point  $x\in X$ we have $\overline{\clM}_{X,x}\simeq \bbN^r$ for some integer $r$,  where  $\overline{\clM}_{X,x}$ denotes the stalk in the \'etale topology.    A morphism between  free monoids $\phi: P_1\to P_2$ is called {\em simple} if $P_1$ and $P_2$ have the same rank, and for every irreducible element of $p_1\in P_1$ there exists a unique element $p_2\in P_2$ and an integer $b$ such that $b\cdot p_2=\phi(p_1)$. A morphism of locally free log structures is called simple if it induces simple morphisms on the stalks. 

Let $D$ be a reduced normal crossing divisor on a scheme $X$.\label{XD-page} According to \cite{KatoLog}, there is a locally free  log structure canonically associated to $D$ in the following way. 
 Let $U:= X\setminus D$ and let $i:U\hookrightarrow X$ be  the inclusion. Then
\begin{equation*}
\clM_D:=i_*(\clO_U^*)\cap \clO_X^*\to \clO_X
\end{equation*}
defines a locally free log structure over $X$. 
Let $x$ be a geometric point of $X$.  The induced morphism 
$$
\overline{\clM}_{D,x}\to \clO_{X,x}
$$
is of the form $\bbN^r\to\clO_{X,x}$ for some integer $r$. In other words,  every irreducible element of the monoid $\overline{\clM}_{D,x}$  corresponds to  an irreducible component of the pullback of $D$ to $\Spec{\clO_{X,x}}$. Roughly speaking, \'etale locally  a normal crossing divisor becomes a {\em simple} normal crossing divisor, namely it is a union of smooth irreducible components. This construction generalizes to stacks.


\subsubsection*{The construction of Matsuki-Olsson}
Let $X$ be a smooth variety and let $D=\cup_{i\in I}D_i\subset X$ be an effective  Cartier divisor with normal crossing support. Let $\{r_i\}_{i\in I}$ be a collection of positive integers.
By \cite{MO},  there exists a smooth  Deligne-Mumford stack $\mathcal{X}$ with a normal crossing divisor
$\mathcal{D}=\cup_{i\in I}\mathcal{D}_{i}\subset \mathcal{X}$ satisfying the following properties:
\label{MO-construction}
\begin{enumerate}
\item The smooth variety $X$ is the coarse moduli space of $\mathcal{X}$.
\item The canonical map  $\pi: \mathcal{X}\to X$ is quasi-finite and flat, and is an isomorphism over $X\setminus D$.
\item $\pi^{*}\mathcal{O}_{X}(-D_i)=\mathcal{O}_{\mathcal{X}}(-r_i\mathcal{D}_{i})$.
\end{enumerate}

Such a stack is defined as a category fibered in groupoids as follows. 
Objects over an $X$-scheme $f:T\to X$ are simple morphisms of log structures
$\phi:f^*\clM_D\to \clM$ such that for any geometric point $t\in T$ with image $x=f(t)\in X$, the induced  morphism  on the stalks of the ghost sheaves is of the following form:
\begin{eqnarray}\label{MO-diagram}
\xymatrix{
\overline{\clM}_{D,\overline{x}}\ar[d]_{\wr}\ar[r]&  \overline{\clM}_{\overline{t}}\ar[d]^{\wr}\\
\bigoplus_{D_i}\bbN\ar[r]^{\oplus(\times r_i)} & \bigoplus_{D_i}\bbN
}
\end{eqnarray}

According to \cite{MO}, if locally $X=\text{Spec}(k[x_1,\cdots,x_n])$ and locally the divisor $D_i=Z(x_i)$ for $1\leq i\leq m$, then $\mathcal{X}$
is canonically isomorphic to the quotient stack
$$[\text{Spec}(k[y_1,\cdots,y_n])/\mu_{r_{1}}\times\cdots\times\mu_{r_{m}}],$$
where $k[y_1,\cdots,y_n]$ is  a $k[x_1,\cdots,x_n]$-algebra via 
$$x_i\mapsto \begin{cases}
y_i^{r_{i}}, & i\leq m \\
y_i, & i> m, 
\end{cases}$$
and the action of $\mu_{r_{1}}\times\cdots\times\mu_{r_{m}}$ is given by 
$$(u_1,\cdots,u_m)\cdot y_i=\begin{cases}
u_iy_i, & i\leq m ,\\
y_i, & i> m.\end{cases}$$

We compare this construction with the root construction. For a smooth scheme $X$, an effective Cartier divisor $D\subset X$, and a positive integer $r$, there exists (see \cite{AGV1}, \cite{Ca}) a smooth Deligne-Mumford stack $X_{(D,r)}$  satisfying the following properties:
\begin{enumerate}
\item  The preimage of $D$ is an infinitesimal neighborhood 
of the $\mu_{r}$-gerbe $\mathcal{D}$ over $D$.
\item There is a canonical map $\pi: X_{(D,r)}\to X$ which is an isomorphism over 
$X\setminus D$.  Every point in $X_{(D,r)}$ lying over $D$ has stabilizer $\mu_{r}$.
\end{enumerate}
This is the $r$-th root construction of $X$ with respect to the divisor $D$ and $r$. 

Let $\mathbb{D}:=(D_1,\cdots,D_n)$ be an $n$-tuple of Cartier divisors and 
$\vec{r}=(r_1,\cdots,r_n)$ be an $n$-tuple of positive integers. Let $X_{(\mathbb{D},\vec{r})}$
be the stack obtained by iterating the root constructions over $X$ and the sequence of divisors. One can see that if the divisor $D=\cup_i D_i$ has simple normal crossing, then $\mathcal{X}\simeq  X_{(\mathbb{D},\vec{r})}$. However if components of $D$ have self-intersections, then along such self-intersections $\clX$ has more automorphisms than $X_{(\mathbb{D},\vec{r})}$. 

\subsubsection{The stack of  twisted curves}\label{subsection:tw-curves-stack}
In \cite{OLogCurv}  the stack of twisted curves $\MM_{g,n}^{tw}$ is constructed using logarithmic geometry. The stack  $\MM_{g,n}^{tw}$ is a smooth Artin stack which has a natural map to the stack of prestable curves  $\MM_{g,n}$ introduced in \cite{Beh97GW}. Such a map is defined by  sending a marked twisted curve $(\clC,\{\sigma_i\})$ to its coarse moduli space with  marked points induced by the $\sigma_i$.  

The notion of {\em log twisted curve} is introduced in \cite{OLogCurv}.
\begin{defn}[\cite{OLogCurv}, Definition 1.7]\label{log_tw_curve_def}
An {\em $n$-pointed log twisted curve} over  a scheme $S$ is a  collection of data
\begin{equation*}
  (C/S,\{\sigma_i,a_i\}, l:\clM_S\to \clM'_S),
\end{equation*}
where $C/S$ is an n-pointed prestable curve, $\sigma_i:S\to C$ are sections (marked points), $a_i$, $i=1,..,n$ are integer-valued locally constant functions on $S$ such that for each $s\in S$ the integer $a_i(s)$ is positive and invertible in the residue field  $k(s)$, and $l:\clM_S\hookrightarrow \clM'_S$ is a  simple morphism of log structures over $S$, where  $\clM_S$ is the canonical log structure associated to $C/S$.
\end{defn}

Log twisted curves turn out to be equivalent to usual twisted curves.

\begin{thm}[\cite{OLogCurv}, Theorem 1.9] For any scheme $S$, there is a natural equivalence of groupoids between the groupoid of $n$-pointed twisted curves over $S$ and the groupoid of $n$-pointed log twisted curves over $S$. Moreover, the equivalence is compatible with base change $S'\to S$. 
\end{thm}

Using this equivalence $\MM_{g,n}^{tw}$ can be seen as the stack over $\bbZ$
which to any $S$ associates the groupoid of $n$-marked  genus $g$  log  twisted curves $(\clC/S, \{\sigma_i, a_i \}, l:\clM_S\hookrightarrow \clM_S')$. For a $n$-tuple of integer numbers $\vec{b}=(b_1,..,b_n)$, let  $\MM_{g,n}^{tw}(\vec{b})$ be the substack of  $\MM_{g,n}^{tw}$ classifying log twisted curves with 
$a_i=b_i$ for all $i$. There is a decomposition in open and closed components
\begin{eqnarray}\label{MM-tw-decomp}
\MM_{g,n}^{tw}\simeq \coprod_{\vec{b}}\MM_{g,n}^{tw}(\vec{b}).
\end{eqnarray}
All the components $\MM_{g,n}^{tw}(\vec{b})$ are isomorphic with each other. The boundary of $\MM_{g,n}(\vec{b})$ is a normal crossing divisor $D$. Then  there is an associated log structure, that we denote by  $\clM_{D}$. For any $\vec{b}$, $\MM_{g,n}^{tw}(\vec{b})$ is the stack over $\MM_{g,n}$ whose fiber over any $f:T\to \MM_{g,n}$ is  the groupoid of 
simple extensions of log structures $f^*\clM_{D}\hookrightarrow \clM_T$ such that for any geometric point $t\in T$ with $f(t)=x$, $\mbox{Coker}(\overline{\clM}^{gp}_{D,x}\to \overline{\clM}^{gp}_{T,t})$ is invertible in $k(t)$.
From the construction just described, we see that  $\MM^{tw}_{g,n}(\vec{b})$ and $\MM_{g,n}$ are locally isomorphic outside of the boundary locus, while over the  locus of singular curves   $\MM^{tw}_{g,n}(\vec{b})$ acquires more automorphisms due to twisted nodes.

Given a log twisted curve $(C/S, \{\sigma_i,a_i\}, l:\clM_S\to \clM'_S)$,  the corresponding  twisted curve $\clC/S$ can be reconstructed as follows. It is the  category fibered in groupoids  whose fiber 
 over any $h:T\to S$ is the groupoid of data consisting of a morphism $s:T\to C$ over $h$ together with a commutative diagram of locally free log structures on $T$:
\begin{eqnarray}\label{from-log-tw-to-tw-diag}
\xymatrix{
h^*\clM_S\ar[r]^l\ar[d] & h^*\clM'_S\ar[d]^{\tau}\\
s^*\clM_C\ar[r]^{k} & \clM'_C,
}
\end{eqnarray}
where 
\begin{enumerate}\label{from-log-tw-to-tw-properties}
\item \label{from-log-tw-to-tw-properties-1}
the morphism $k$ is simple and for any geometric point $t$ of $T$, the map $\overline{\clM}'_{S,t}\to~\overline{\clM}'_{C,t}$ is either an isomorphism, or of the form $\bbN^r\to \bbN^{r+1}$ mapping $e_i$ to $e_i$ for $i < r$ and $e_r$ to either $e_r$ or $e_r+e_{r+1}$, and  
\item  \label{from-log-tw-to-tw-properties-2}
for every $1\leq i\leq n$ and geometric point $t$ of $T$ with image $s=s(t)$ in $\sigma_i(S)\subset C$, the group $$Coker(\overline{\clM}'^{gp}_{S,t}\oplus   \overline{\clM}^{gp}_{C,t}\to  \overline{\clM}'^{gp}_{C,t})$$ is a cyclic group of order $a_i$. 
\end{enumerate}


\section{Moduli of twisted  stable maps to root  gerbes}\label{moduli_maps_root_gerbe}

Let $X$ be a smooth projective variety over $\bbC$, $\clL$ a line bundle over $X$, and $r\geq 1$ an integer. The purpose of this Section is to study the structure of the moduli stack $\clK_{0,n}(\clG, \beta)$ of genus $0$ twisted stable maps to a root gerbe $\clG:=\sqrt[r]{\clL/X}$. More precisely, we study components\footnote{The term ``component'' here means a union of connected components.} $\clK_{0,n}(\clG, \beta)^{\vec{g}}$ of $\clK_{0,n}(\clG, \beta)$ indexed by what we call {\em $\beta$-admissible vectors} (Definition \ref{adm_vector-def}). The main results of this Section, Theorems \ref{mod_space_to_gerbe_is_gerbe} and \ref{mod_space_to_root_gerbe_is_root_gerbe}, exhibit the structure of $\clK_{0,n}(\clG, \beta)^{\vec{g}}$ over the moduli stack\footnote{We always assume that $\overline{M}_{0,n}(X, \beta)$ is non-empty.} $\overline{M}_{0,n}(X, \beta)$ of stable maps to $X$. 

\subsection{Components of moduli stack}
We begin with some useful Lemmas.
\begin{lem}
\label{gerbe_lift}
 Let $G$ be a finite group, and let $\clG\to X$ be a $G$-banded gerbe. To give a lift $\clC\to\clG$ of a map to the coarse moduli space $C\to X$ is equivalent to give a map $\clC\to \clG\times_X C$.
\end{lem}
\begin{pf}
Consider  the following diagram
\begin{eqnarray}
\label{map-gerby-curve-diag}
\xymatrix{
\clC \ar[dr]_{\pi}\ar[r]_-{\overline{f}}\ar@/^1pc/[rr]^{\tilde{f}}& \clG\times_X C\ar[r] \ar[d] \ar@{}[rd]|{\square}& \clG\ar[d]\\
 &  C\ar[r]^{f} & X.
}
\end{eqnarray}
To give a lift  $\tilde{f}$ of $f$  means to give the outer square in diagram
(\ref{map-gerby-curve-diag}), namely the pair $(\tilde{f},\pi)$.
Due to the universal property of the fiber product this   is equivalent to give a  map $\overline{f}: \clC\to \clG\times_X C$. Moreover  $\overline{f}$  is representable if and only if $\tilde{f}$ is.
\end{pf}

\begin{lem}[c.f. \cite{Ca}]\label{rep_map_to_root}
Let $(\clC, \{\sigma_i\})$ be an $n$-pointed smooth twisted curve with stack points $\sigma_i, 1\leq i\leq n$. Let $\mu_{r_i}$ be the isotropy group of the stack point $\sigma_i$. Denote by $\pi:(\clC, \{\sigma_i\})\to~(C,\{p_i\})$ the coarse curve. Let  $\tilde{f}:\clC\to\sqrt[r]{\clL/X}$ be a morphism and $f:C\to~X$ its induced map between coarse moduli spaces. Suppose $\tilde{f}$ is given by a line bundle $M= \pi^*L\otimes \bigotimes_{i=1}^n \clT_i^{m_i}$ over $\clC$ (with $0\leq m_i< r_i$) and an isomorphism $\psi: M^{\otimes r}\simeq \pi^*f^*\clL$. Then $\tilde{f}$ is representable if and only if for $1\leq i\leq n$, we have $r_i|r$ and $m_i$ and $r_i$ are co-prime. 
\end{lem}
\begin{pf}
By \cite{AV}, Lemma 4.4.3, it suffices to study the homomorphism 
\begin{equation}\label{hom_between_aut}
Aut(\sigma_i)\to Aut(\tilde{f}(\sigma_i)),
\end{equation}
 induced by $\tilde{f}$ on stack points. Here by $\sigma_i$ we mean a morphism 
$\tilde{h}_i: \Spec{K}\to \clC$ from  an algebraically closed field  $K$ to $\clC$ with image in   the special locus. 
 By the root construction description of $\clC$ (see e.g. \cite{Ca}, Example 2.7 and \cite{AGV1}, Section 4.2), the stack point $\sigma_i$ is equivalent to the data $(h_i, M_i, t_i,\phi_i)$, where $h_i:\Spec{K}\to C$ with image $p_i$, $M_i$ is a line bundle over $\Spec{K}$,  $\phi_i:M_i^{\otimes r_i}\stackrel{\sim}{\to} h_i^*\clO(p_i)$, $t_i$ is a section of $M_i$ such that $\phi_i(t_i^{r_i})= h_i^* s_i$, hence $t_i=0$. The image $\tilde{f}(\sigma_i)$ is given by $\tilde{h}_i^*M$ and $\tilde{h}_i^*\psi: \tilde{h}_i^*M^{\otimes r}\simeq \tilde{h}_i^*\pi^*f^*\clL$. 
Note that $\tilde{h}_i^*\clT_i$ is naturally isomorphic to $M_i$. An automorphism  $\epsilon\in Aut(\sigma_i)\simeq \mu_{r_i}$ is mapped to $\epsilon^{m_i}\in Aut(\tilde{f}(\sigma_i)))\simeq \mu_r$ since   $M=\pi^*L\otimes \bigotimes_{i=1}^n \clT_i^{m_i}$ .
This homomorphism is injective if and only if $r_i|r$ and $m_i$ and $r_i$ are coprime.
\end{pf}

\subsubsection{Admissible vectors}
The inertia stack $I\clG$ admits a decomposition $$I\clG=\cup_{g\in \mu_r} \clG_g$$ indexed by elements of $\mu_r$. An object of $\clG_g$ over $h:T\to X$ is a collection $((M, \phi), g)$ where $(M, \phi)$ is an object of $\clG$ over $T$ (i.e. $M$ is a line bundle over $T$ and $\phi: M^{\otimes r}\to h^*\mathcal{L}$ is an isomorphism) and $g$ is an automorphism of $(M, \phi)$ defined by multiplying\footnote{The identification of $\mu_r$ with the group of $r$-th roots of $1\in \bbC^*$ allows us to identify $g\in \mu_r$ with complex numbers. We use this to make sense of the multiplication. In what follows we use this identification without explicit reference.} fibers of $M$ by $g$. 

\begin{defn}\label{adm_vector-def}
Let $\bar{I}(\clG)_g\subset \bar{I}(\clG)$ be the image of $\clG_g$ under the natural map $I\clG\to \bar{I}(\clG)$. Let $\vec{g}:=(g_1,...,g_n)\in \mu_r^{\times n}$ be a vector of elements of $\mu_r$. Set 
\begin{eqnarray}
\clK_{0,n}(\clG,\beta)^{\vec{g}}:=\cap_{i=1}^n ev_i^{-1}( \overline{I}(\clG)_{g_i})
\end{eqnarray}
The vector $\vec{g}$ is called {\em $\beta$-admissible} if $\clK_{0,n}(\clG,\beta)^{\vec{g}}$ is nonempty. 
\end{defn}
\begin{rmk}
Note that the definition of $\beta$-admissible vectors depends on a choice of the class $\beta$.
\end{rmk}

Let $[\tilde{f}: (\clC, \{\sigma_i\})\to \clG]\in \clK_{0,n}(\clG, \beta)^{\vec{g}}(\bbC)$. By definition the morphism $\tilde{f}|_{\sigma_i}:B\mu_{r_i}\simeq \sigma_i\to \clG$ is equivalent to an injective  homomorphism  $$\mu_{r_i}\hookrightarrow \mu_r, \quad \exp(2\pi \sqrt{-1}/r_i)\mapsto g_i.$$
The argument in the proof of Lemma \ref{rep_map_to_root}, applied to the irreducible component of $\clC$  containing $\sigma_i$, shows that we may write 
\begin{equation}\label{elements_in_adm_vector}
g_i=\exp(2\pi\sqrt{-1}\frac{m_i}{r_i}), \quad \text{with } 0\leq m_i< r_i, \text{ and } (m_i, r_i)=1.
\end{equation}
Furthermore, if $\clL^{1/r}$ is the universal $r$-th root of $\clL$ over $\clG$, then $\tilde{f}|_{\sigma_i}^*\clL^{1/r}$ is the $\mu_{r_i}$-representation on which the standard generator $\exp(2\pi\sqrt{-1}/r_i) \in \mu_{r_i}$ acts by multiplication by $\exp(2\pi\sqrt{-1}m_i/r_i)$. In other words 
\begin{equation}\label{age_of_universal_root_line_bundle}
\text{age}_{\sigma_i}(\tilde{f}^*\clL^{1/r})=\frac{m_i}{r_i}. 
\end{equation}

\begin{lem}
\label{adm_vect_lemm}
Suppose $\vec{g}=(g_1,...,g_n)\in \mu_r^{\times n}$ is a $\beta$-admissible vector. Then 
\begin{equation}\label{phase_shift_eq}
\prod_{i=1}^n g_i=\exp\left(\frac{2\pi \sqrt{-1}}{r}\int_\beta c_1(\clL)\right).
\end{equation}
\end{lem}
\begin{pf}
Let $[\tilde{f}:(\clC, \sigma_1,...,\sigma_n) \to \sqrt[r]{\clL/X}]\in \clK_{0,n}(\clG,\beta)^{\vec{g}}(\bbC)$ be a twisted stable map. Let $\clL^{1/r}$ be the universal $r$-th root of $\clL$ over $\clG=\sqrt[r]{\clL/X}$. By Riemann-Roch for twisted curves (see e.g. \cite{AGV1}, Theorem 7.2.1),
\begin{eqnarray}\label{RR_tw_eq}
\chi(\tilde{f}^*\clL^{1/r})=1+ \text{deg}\, \tilde{f}^*\clL^{1/r}-\sum_{i=1}^{n} \text{age}_{\sigma_i}(\tilde{f}^*\clL^{1/r}),
\end{eqnarray}
which is an integer. Clearly $$\text{deg}\, \tilde{f}^*\clL^{1/r}=\frac{1}{r}\int_{\beta}c_1(\clL).$$
By (\ref{elements_in_adm_vector}) and (\ref{age_of_universal_root_line_bundle}) we have $$g_i=\exp(2\pi \sqrt{-1}\ \text{age}_{\sigma_i}(\tilde{f}^*\clL^{1/r})).$$ The result follows.
\end{pf}


\begin{prop}
\label{unique_geom_lift}
Let $\clG=\sqrt[r]{\clL/X}\to X$ be a root gerbe. Let $[f:(C, p_1,...,p_n) \to X]$ be an object of $\overline{M}_{0,n}(X,\beta)(\bbC)$.  Then for a vector $\vec{g}=(g_1,...,g_n)\in \mu_r^{\times n}$ satisfying (\ref{phase_shift_eq}) there exists, up to isomorphisms, a unique twisted stable map $\tilde{f}:(\clC, \sigma_1,..., \sigma_n)\to \clG$ in $\clK_{0,n}(\clG,\beta)^{\vec{g}}$ lifting $f$.
\end{prop}

\begin{pf}
We first assume that $C$ is smooth. Associate to $\vec{g}$ the numbers $r_i$ and $m_i, 1\leq i\leq n$ as in (\ref{elements_in_adm_vector}). Let $(\clC, \sigma_1,...,\sigma_n )$ be the smooth twisted curve obtained by applying the $r_i$-th root construction to the divisor $p_i\subset C$ for $1\leq i\leq n$. Denote by $\clT_i, 1\leq i\leq n$ the tautological sheaves and by $\pi: \clC\to C$ the natural map. By (\ref{phase_shift_eq}) we have $$\frac{1}{r}\text{deg}\, \pi^*f^*\clL-\sum_{1\leq i\leq n}\frac{m_i}{r_i}\in \mathbb{Z}.$$ Pick $L\in \text{Pic}(C)$ such that $\text{deg}\, L=\frac{1}{r}\text{deg}\, \pi^*f^*\clL-\sum_{1\leq i\leq n}\frac{m_i}{r_i}$. Set $$M:=\pi^*L\otimes \bigotimes_{1\leq i\leq n} \clT_i^{m_i}.$$ Then $\text{deg}\, M^{\otimes r}=\text{deg}\, \pi^*f^*\clL$, so there exists an isomorphism $M^{\otimes r}\simeq \pi^*f^* \clL$, which defines a map $\tilde{f}: \clC\to \clG$. By construction $\tilde{f}$ is a lifting of $f$. By Lemma \ref{rep_map_to_root} $\tilde{f}$ is representable. Also, $\tilde{f}$ is unique up to isomorphisms since the line bundle $L$ on $C\simeq \mathbb{P}^1$ is determined up to isomorphisms by its degree. This proves the Proposition in case $C$ is smooth. 

We treat the general case by induction on the number of irreducible components of $C$. The case of one irreducible component is proven above. We now establish the induction step. Let $C_1\subset C$ be an irreducible component containing only one node $x$, and $C_2:=\overline{C\setminus C_1}$. In other words $C_1$ is an irreducible component meeting the rest of the curve $C_2$ at the node $x$. Let $T\subset [n]:=\{1,2,...,n\}$ be the marked points that are contained in $C_1$. Restrictions of $f$ yield two stable maps $$f_1: (C_1, \{p_i| i\in T\}\cup \{x\})\to X, \quad f_2: (C_2, \{p_i|i\in T^C\}\cup \{x\})\to X.$$ Here $T^C:=[n]\setminus T$. Set $\beta_1:=f_{1*}[C_1]$ and define $m_x, r_x\in \mathbb{Z}$ by $$\frac{m_x}{r_x}:=\langle-\sum_{i\in T}\frac{m_i}{r_i}+\frac{1}{r}\int_{\beta_1}c_1(\clL) \rangle, \quad (m_x, r_x)=1.$$ Here $\langle - \rangle$ denotes the fractional part. Note that $r_x|r$ since $r_i|r$ for $i\in T$. Therefore we may define $g_x:=\exp(2\pi\sqrt{-1}\frac{m_x}{r_x})\in \mu_r$. 

By the smooth case there exists a lifting $$\tilde{f}_1: (\clC_1, \{\sigma_i| i\in T\}\cup \{\sigma_x\})\to \clG$$ of $f_1$ associated to the collection $\{g_i|i\in T\}\cup \{ g_x\}$. By induction there exists a lifting $$\tilde{f}_2: (\clC_2, \{\sigma_i| i\in T^C\}\cup \{\sigma'_x\})\to \clG$$ of $f_2$ associated to the collection $\{g_i| i\in T^C\}\cup \{g_x^{-1}\}$. By our choices of the actions of isotropy groups at $\sigma_x$ and $\sigma'_x$ we see that $\clC_1$ and $\clC_2$ glue along $\sigma_x, \sigma'_x$ to form a balanced twisted curve $\clC$, and the morphisms $\tilde{f}_1, \tilde{f}_2$ define a morphism $\tilde{f}: (\clC, \sigma_1,..., \sigma_n)\to \clG$ which is representable and is a lifting of $f$. 

By restricting to irreducible components, we see that uniqueness of lifting follows from uniquess of lifting in the smooth case. This completes the proof.    
\end{pf}

\begin{rmk}
Lemma \ref{adm_vect_lemm} and Proposition \ref{unique_geom_lift} combined show that $\beta$-admissible vectors (for a fixed class $\beta$) are completely characterized by the condition (\ref{phase_shift_eq}). This condition can be viewed as a generalization of the monodromy condition which is required to hold for genus $0$ twisted stable maps to $B\mu_r$. 
\end{rmk}


In what follows we prove that each open-and-closed component $\clK_{0,n}(\clG,\beta)^{\vec{g}}$ of  $\clK_{0,n}(\clG,\beta)$ is a gerbe over a base stack that can be constructed over $\overline{M}_{0,n}(X,\beta)$ by using logarithmic geometry. This base stack is isomorphic to $\overline{M}_{0,n}(X,\beta)$ over the  (possibly empty) locus corresponding to twisted stable maps with smooth domain curve. Along the boundary it has more automorphisms, corresponding to the  fact that singular twisted curves carry additional automorphisms associated to the stacky nodes. In order to describe the gerbe structure of $\clK_{0,n}(\clG,\beta)^{\vec{g}}$ and its base stack we will need an auxiliary   stack    parametrizing weighted  prestable curves.  

\subsection{The stacks $\mathfrak{M}_{g,n,\beta}$ and  $\mathfrak{M}^{tw}_{g,n,\beta}$ }\label{frak-M-weight-def}
We describe  a stack  introduced in \cite{Cost03}, Section 2. It parametrizes  weighted prestable curves with a stability condition. Denote by $\MM_{g,n}$ the stack  of genus $g$ prestable curves with $n$ marked points.

Let $g, n\in \bbZ_{\geq 0}$ and $\beta\in H_2^+(X,\bbZ)$. A triple $(g,n,\beta)$ is called {\em stable} if either $\beta\neq 0$ or $\beta=0$ and $2g-2+ n>0$. For a triple $(g,n,\beta)$ the stack $\mathfrak{M}_{g,n,\beta}$ over $\mathfrak{M}_{g,n}$ is defined inductively in the following way:
\begin{enumerate}
\item If $(g,n,\beta)$  is unstable, then $\mathfrak{M}_{g,n,\beta}$  is empty.
\item If  $(g,n,\beta)$  is stable, an object of $\mathfrak{M}_{g,n,\beta}$ over $T$ is 
\begin{enumerate}
\item  an object $(C,\{s_i\})$  of  $\mathfrak{M}_{g,n}(T)$, namely a genus $g$    prestable curve over  $T$  with  $n$ sections  in the smooth locus;
\item  a constructible function $f: C_{gen}\to H_2^+(X,\bbZ)$, where  $C_{gen}\to T$  is   the complement of the nodes and the sections in $C$.  The function $f$  must be locally constant on the geometric fibers of $C_{gen}\to T$.
\end{enumerate}
\item if $T^0\subset T$ is the open subscheme parametrizing nonsingular curves $C^0\to T^0$ , then $f:C^0_{gen}\to H_2^+(X,\bbZ)$ must be constant
with value $\beta$.

\item $f$ has to satisfy two kinds of  gluing conditions  along the boundary of  $\mathfrak{M}_{g,n}$:
\begin{enumerate}
\item Suppose that  there is a decomposition $g=g'+g''$ and 
$[n]=\{1,..,n\}=T\coprod T^c$, with $\abs{T}=n'$ and $\abs{T^c}=n''$    and a map $S\to T$ such that the composite map $S\to \mathfrak{M}_{g,n}$ factors into
\begin{equation*}
S\to \mathfrak{M}_{g',T\coprod\{s'\}}\times  \mathfrak{M}_{g'',T^C\coprod\{s''\}}\to \mathfrak{M}_{g,n},
\end{equation*}
where the second map is obtained by gluing the marked sections $s'$ and $s''$. 
Let  $C_V' \in\mathfrak{M}_{g',n'}(S)$   and  $C_V''\in \mathfrak{M}_{g^{''},n''}(S)$ be the associated families of curves. We require that  the pulled back constructible functions  $f':C_V'\to H_2^+(X,\bbZ)$ and  $f'':C_V''\to H_2^+(X,\bbZ)$  define a morphism
\begin{equation*}
S\to \coprod_{\beta'+\beta''=\beta}\mathfrak{M}_{g',T,\beta'}\times\mathfrak{M}_{g'',T^c,\beta''}. 
\end{equation*} 
\item  Suppose that there is a map $S\to T$ such that the composite map $S\to \mathfrak{M}_{g,n}$ factors into
\begin{equation*}
S\to \mathfrak{M}_{g-1,n \coprod\{s',s''\}}\to  \mathfrak{M}_{g,n}.
\end{equation*}
Then the associated  genus $g-1$ family of curves $C_S\to S$ with the pulled back constructible function $f:C_{Sgen}\to H_2^+(X, \bbZ)$ has to define a morphism
\begin{equation*}
S\to  \mathfrak{M}_{g-1,n\coprod\{s',s''\},\beta}.
\end{equation*}
\end{enumerate}
\end{enumerate}
The constructible function $f: C_{gen}\to H_2^+(X, \bbZ)$ will be called the {\em weight}. An object $(C, f)$ of $\MM_{g,n, \beta}$ is called a $H_2^+(X, \bbZ)$-weighted prestable curve, and its {\em total weight} is by definition $\beta$. 

Note that in this definition it is important  that $H_2^+(X,\bbZ)$ is a semigroup with indecomposable zero and any of its elements has a finite number of decompositions.   

We quote some properties of the stack $\MM_{g,n, \beta}$.
\begin{prop}[\cite{Cost03}, Proposition 2.0.2]
The map  $\mathfrak{M}_{g,n\beta}\to \mathfrak{M}_{g,n}$ defined by forgetting the weights is \'etale, and relatively a scheme of finite type. Therefore the stack  $\mathfrak{M}_{g,n,\beta}$ is a smooth algebraic stack. 
\end{prop}
\begin{prop}[\cite{Cost03}, Proposition 2.1.1]\label{universal_family}
The natural morphism $\MM_{g,n+1, \beta}\to \MM_{g,n,\beta}$ defined by forgetting the $(n+1)$-st marked point is the universal family over $\MM_{g,n,\beta}$.  
\end{prop}

\subsubsection{Boundary of $\mathfrak{M}_{0,n,\beta}$}
\label{bdry-frak-M-weigh-sect}
Boundary divisors of $\MM_{0,n}$ are indexed by subsets $T$ of $[n]:=\{1,2,...,n\}$. Each $T$ corresponds to the boundary divisor $D^T$ which parametrizes curves $C=C_1\cup C_2$ meeting at a node such that marked points indexed by $T$ are contained in $C_1$ and other marked points are contained in $C_2$. 

Boundary divisors in $\MM_{0,n,\beta}$ can be similarly described.
\begin{defn}\label{boundary_of_M0nbeta}
Given $(T, \beta')$, where $T$ is a not necessarily proper subset of $[n]$ and $\beta'\in H_2^+(X,\bbZ)$ such that $\beta'\leq \beta$, define $D^T_{\beta'}\subset \mathfrak{M}_{0,n,\beta}$ to be the divisor which parametrizes curves $C=C_1\cup C_2$  meeting at a node such that 
\begin{enumerate}
\item
marked points indexed by $T$ are contained in $C_1$,  marked points indexed by $T^C:= [n]\setminus T$ are contained in $C_2$;  
 \item
 $f|_{C_1}=\beta'$ and $f|_{C_2}=\beta-\beta'$.
 \end{enumerate} 
\end{defn}
 
 Note that $D^T_{\beta'}=D^{T^C}_{\beta-\beta'}$. Let $h_n:\MM_{0,n+1,\beta}\to \MM_{0,n,\beta}$ be the map that forgets the $(n+1)$-th marked point. Then 
\begin{equation}\label{preimage-of-divisor-under-univ-curve}
h_n^{-1}D^T_{\beta'}=D^T_{\beta'}\cup D^{T\cup\{n+1\}}_{\beta'}.
\end{equation}

\begin{lem}
The boundary divisors $D^T_{\beta'}$ are normal crossing divisors.
\end{lem}
\begin{pf}
Consider the natural map $l_n:\MM_{0,n,\beta}\to \MM_{0,n}$ that forgets the weights. The following relation holds
\begin{equation*}
l_n^{-1} D^T=\begin{cases}\cup_{0 < \beta'\leq \beta} D^T_{\beta'}\quad  \text{if } \abs{T}<2;\\
\cup_{0\leq \beta'\leq \beta} D^T_{\beta'}\quad  \text{if } 2\leq \abs{T}\leq  n- 2;\\
 \cup_{0\leq \beta' < \beta} D^T_{\beta'}\quad \text{if }  \abs{T} > n- 2 .\end{cases}
\end{equation*}
Since $D^T$ are normal crossing divisors in $\MM_{0,n}$ and the morphism $l_n$ is \'etale, the result follows.
\end{pf}

\label{log_str_MM_0nbeta_page}
\subsubsection{Log structure on $\MM_{0,n,\beta}$.}
Let  $\mathcal{I}_D$ \label{calI_D_def_intext} be  the set of pairs $(T, \beta')$ with $T\subset [n], n\notin T$, such that one of the following holds:
\begin{enumerate}\label{set_I_sub_D}
\item
$0< \beta' < \beta$;
\item
$\beta'=0$ and $|T|\geq 2$;
\item
$\beta'=\beta$ and $|T|\leq n-2$.
\end{enumerate}

The union\footnote{We want each divisor $D^T_{\beta'}$ to appear only once in this union. I.e. the pairs  $(T, \beta'), (T^C, \beta-\beta')$ shouldn't both occur. The requirement $n\notin T$ is a way to select only one of them.\label{exclude_double_count_foot}} of boundary divisors  
\begin{equation}\label{union_of_boundary_divisors}
\bigcup_{(T,\beta')\in \mathcal{I}_D  }D^{T}_{\beta'}\subset \MM_{0,n,\beta},
\end{equation}
 is a reduced normal crossing divisor. 
 The divisor (\ref{union_of_boundary_divisors}) defines a locally free log structure over $\MM_{0,n,\beta}$ which we denote by $\clM_D^n$. This follows from a  general construction we described in Section \ref{subsection:log-geom-and-tw-curves} (see page \pageref{XD-page}), following \cite{KatoLog}. 


Let $\MM_{g,n}^{tw}$ be the stack of genus $g$ twisted curves with $n$ marked points  introduced in  Section  \ref{subsection:tw-curves-stack}.  Define $$\MM_{g,n, \beta}^{tw}:=\MM_{g,n, \beta}\times_{\MM_{g,n}}\MM_{g,n}^{tw}.$$
The stack $\MM_{g,n, \beta}^{tw}$ parametrizes $H_2^+(X, \bbZ)$-weighted genus $g$ twisted curves with $n$ marked points and total weight $\beta$.

\subsection{The stack  $\mathfrak{Y}_{0,n,\beta}^{\vec{g}}$}\label{section:stack_YY}
We will define a stack over $\MM_{0,n,\beta}$ using the log structure $\clM_D^n$ and some additional data coming from a $\beta$-admissible vector $\vec{g}$, following \cite{MO}. 

Fix a $\beta$-admissible vector $\vec{g}=(g_1,..., g_n)\in \mu_r^{\times n}$. We define a collection of triples of integers $\{(\rho_i, r_i, m_i)|1\leq i\leq n\}$ as follows. Each $g_i, 1\leq i\leq n$ may be identified with a root of unity $$g_i=\exp(2\pi \sqrt{-1} \theta_i), \quad  \text{where } \theta_i\in \mathbb{Q}\cap [0, 1),$$
which defines the rational numbers $\theta_i, 1\leq i\leq n$. The characterizing relation of $\beta$-admissible vectors (\ref{phase_shift_eq}) reads $$\prod_{i=1}^n \exp(2\pi \sqrt{-1}\theta_i)= \exp(\frac{2\pi \sqrt{-1}}{r} \int_\beta c_1(\mathcal{L})).$$
For $1\leq i\leq n$, define 
\begin{equation}\label{the_triple_of_integers}
\rho_i:= r\theta_i, \quad r_i:=\frac{r}{gcd (r, \rho_i)}, \quad m_i:=\frac{\rho_i}{gcd (r, \rho_i)}.
\end{equation}

Let $(T, \beta')$ be an index of the boundary divisors of $\MM_{0,n,\beta}$ as in Definition \ref{boundary_of_M0nbeta}. Define 
\begin{equation}\label{triple_of_integers_for_node}
\theta_{T, \beta'}:=\langle\frac{1}{r}\int_{\beta'} c_1(\mathcal{L})-\sum_{i\in T} \theta_i\rangle, \quad r_{T, \beta'}:=\frac{r}{gcd(r, r\theta_{T, \beta'})}, \quad m_{T, \beta'}:=\frac{r\theta_{T, \beta'}}{gcd(r, r\theta_{T, \beta'})}.  
\end{equation}
Here $\langle -\rangle$ again denotes the fractional part. This definition makes sense since $\int_{\beta'}c_1(\mathcal{L})-\sum_{i\in T} r\theta_i$ is an integer.

\begin{defn}
Let $\YY_{0,n,\beta}^{\vec{g}}$ be the stack obtained by applying the construction of \cite{MO}, Theorem 4.1 recalled in Section \ref{subsection:log-geom-and-tw-curves}  to the stack $\MM_{0,n,\beta}$, the normal-crossing divisor (\ref{union_of_boundary_divisors}), and the collection of positive integers $\{r_{T, \beta'}| (T, \beta')\in \mathcal{I}_D\}$. 

Let $\YY_{0,n+1, \beta}^{\vec{g}\cup \{1\}}$ be the stack obtained in the same way from the stack $\MM_{0, n+1, \beta}$ and the $\beta$-admissible vector $\vec{g}\cup\{1\}:= (g_1,...,g_n, 1)\in \mu_r^{\times n+1}$. 
\end{defn}
As in the proof of \cite{MO}, Theorem 4.1, the stack $\YY_{0,n,\beta}^{\vec{g}}$ is defined as a category fibered in groupoids whose objects over a $\MM_{0,n, \beta}$-scheme  $f: S\to\MM_{0,n, \beta}$ are  simple morphisms of log structures 
$$f^*\clM_D^n\to\clM_S$$ such that for every geometric point $\overline{s}\to S$ of $S$ with $x=f(\overline{s})$ there is a commutative diagram
\begin{eqnarray}
\label{Mats_Olss_diag}
\xymatrix{
\bigoplus_{(T_i, \beta_i)}\bbN\ar[rr]^-{\oplus (\times r_{T_i,\beta_i})}\ar[d]_{\wr} &  & \bigoplus_{(T_i,\beta_i)}\bbN \ar[d]^{\wr}\\
f^{-1}\overline{\clM}_{D,\overline{x}}^n\ar[rr] & &\overline{\clM}_{S,\overline{s}},
}
\end{eqnarray}
where the sum is taken over  pairs $(T_i,\beta_i)\in\mathcal{I}_D$ labelling  the irreducible components
of the pullback of $D$ to $\Spec{\clO_{S,\overline{s}}}$. Note that the  $(T_i,\beta_i)$ may have repetitions. 
\begin{rmk}\label{YY_to_MM-flat}
Note that property (2) listed in the construction of Matsuki-Olsson in Section  \ref{MO-construction} implies that the structure morphism $\mathfrak{Y}_{0,n,\beta}^{\vec{g}}\to\MM_{0,n,\beta}$ is quasi-finite and flat.
\end{rmk}

\begin{prop}\label{description_of_stack_YY}
The stack  $\mathfrak{Y}_{0,n,\beta}^{\vec{g}}$  parametrizes $H_2^+(X,\bbZ)$-weighted genus $0$ twisted curves such that 
\begin{itemize}
\item
the $i$-th marked gerbe is banded by $\mu_{r_i}$;
\item
the nodes are stack points of orders $\{r_{T, \beta'}\}$.
\end{itemize}
In particular $\YY_{0,n, \beta}^{\vec{g}}$ is an open substack of $\MM_{0,n, \beta}^{tw}$.
\end{prop}
\begin{pf}
By definition a morphism $\tilde{f}: S\to \mathfrak{Y}_{0,n,\beta}^{\vec{g}}$ consists of the folowing data:
\begin{enumerate}
\item
a morphism $f: S\to\MM_{0,n,\beta}$, corresponding to a weighted curve $(C/S, \{s_i\})$ with marked sections $s_i$; 
\item
a simple morphism of log structures $l: f^*\clM^n_D\hookrightarrow \clM_S$  inducing diagram (\ref{Mats_Olss_diag}),   with coefficients associated to irreducible components of $D$ as in (\ref{triple_of_integers_for_node});
\item a $n$-tuple of integer numbers $\{r_i\}_{i=1}^n$, determined by $\vec{g}$ and $\beta$ according to (\ref{the_triple_of_integers}).
\end{enumerate}
The above data are equivalent by definition to a log twisted curve, see Definition \ref{log_tw_curve_def}. The corresponding twisted curve $\clC$ is determined as in \cite{OLogCurv}, Section 4 (see here page \pageref{from-log-tw-to-tw-diag}).  Consider  a morphism $U\to \clC$. 
Let $\overline{u}\to U$ be a geometric point mapping to the marked point $s_i$ of $C$. Let $t$ be an element in $\clO_{C,\overline{u}}$ locally defining $s_i$. According to \cite{OLogCurv}, Section 4.2, \'etale locally $\clC$ is isomorphic to
\begin{equation}
\left[\Spec{(\clO_{C,\overline{u}}[z]/(z^{r_i}-t))}/\mu_{r_i}\right],
\end{equation}
where $\mu_{r_i}$ acts by multiplication on $z$. Let $\overline{u}\to C$ map to a node. Such a node corresponds to an irreducible component $D^{T_i}_{\beta_i}$ of the boundary divisor $D$ on $\MM_{0,n,\beta}$.
Let us  consider the pullback of   $D^{T_i}_{\beta_i}$ to $\Spec{\clO_{U,\overline{u}}}$.  If the curve $C$ has $k$ nodes of type $(T_i,\beta_i)$, the pulled back divisor will have $k$ irreducible components corresponding to 
 irreducible elements  $e_{i_1},..,e_{i_k}$ of the monoid $f^{-1}\overline{\clM}^n_{D,\overline{u}}\simeq \bbN^r$, where $r$ is some integer number.
 The induced morphism   $l_{\overline{u}}: f^{-1}\overline{\clM}^n_{D,\overline{u}}\to \overline{\clM}_{S,\overline{u}}$ acts on  the submonoid generated by  $e_{i_1},..,e_{i_k}$ as the multiplication by  $r_{T_i,\beta_i}$.
  According to \cite{OLogCurv}, Section 4.3, after choosing an \'etale morphism $C\to \Spec{\clO_{S,\overline{u}}[x,y]/(xy-t)}$, the \'etale local description of $\clC$ is   
\begin{equation}
[\Spec{(\clO_{C,\overline{u}}[z,w]/(zw-t', z^{r_{T_i,\beta_i}}=x,  w^{r_{T_i,\beta_i}}=y  ))/\mu_{r_{T_i,\beta_i}}}],
\end{equation}
where the action  of $\mu_{r_{T_i,\beta_i}}$ is the usual balanced action. 
\end{pf}

\begin{lem}\label{some_map_is_etale}
The natural morphism $\YY_{0,n,\beta}^{\vec{g}}\to \MM_{0,n,\beta}^{tw, \vec{g}}$  is \'etale.
\end{lem}
\begin{pf}
This is immediate from the proof of Proposition \ref{description_of_stack_YY} and \cite{OLogCurv}, Lemma 5.3. 
\end{pf}

By Proposition \ref{description_of_stack_YY} the stack $\YY_{0,n,\beta}^{\vec{g}}$ is the moduli stack of certain weighted twisted curves. Its universal family can be described as follows.

\begin{prop}\label{family_of_twisted_curves-lem}
There exists a natural morphism 
\begin{equation}\label{map_from_YYn+1_to_YYn}
\YY_{0, n+1, \beta}^{\vec{g}\cup\{1\}}\to \YY_{0,n,\beta}^{\vec{g}},
\end{equation}
 which is the universal family of genus $0$ twisted curves with $n$ marked gerbes over $\YY_{0,n,\beta}^{\vec{g}}$. In particular for $1\leq i\leq n$, the $i$-th marked gerbe is banded by the group $\mu_{r_i}$.
\end{prop}

\begin{pf}
 By the construction of $\YY_{0,n,\beta}^{\vec{g}}$, it has a normal crossing divisor 
\begin{equation*}
\clD=\bigcup_{(T,\beta')\in \mathcal{I}_D} \clD_{\beta'}^{T}, 
\end{equation*}
such that $\pi_n^*\mathcal{O}_{\MM_{0,n,\beta}}(-D^T_{\beta'})=\mathcal{O}_{\YY_{0,n,\beta}^{\vec{g}}}(-r_{T, \beta'}\clD^T_{\beta'})$, where $\pi_n: \YY_{0,n,\beta}^{\vec{g}}\to \MM_{0,n,\beta}$ is the natural map. Let $\clM_{\clD}^n$ be the locally free log structure associated to $\clD$. By construction of $\YY_{0,n,\beta}^{\vec{g}}$ there is a universal simple morphism 
\begin{equation}\label{frak-l-universal-simple-morph-eq}
\mathfrak{l}_n: \pi_n^*\clM_D^n\to \clM_{\clD}^n,
\end{equation} 
of  log structures over  $ \YY_{0,n,\beta}^{\vec{g}}$, where $\clM_D^n$ is the log structure defined by the divisor (\ref{union_of_boundary_divisors}).  

Let $h_n: \CC_{0,n,\beta}\to \MM_{0,n,\beta}$ be the universal family of weighted curves over $\MM_{0,n,\beta}$. By Proposition \ref{universal_family}, $\CC_{0,n,\beta}$ can be identified with $\MM_{0,n+1,\beta}$. Let $h_n': \CC'\to \YY_{0,n,\beta}^{\vec{g}}$ be the pull-back of $h_n$ via the natural map $\pi_n: \YY_{0,n,\beta}^{\vec{g}}\to \MM_{0,n,\beta}$, i.e. there is a $2$-cartesian diagram
\begin{equation*}
\xymatrix{
\CC'\ar[r]^{\pi'_{n+1}}\ar[d]_{h'_n}\ar@{}[rd]|{\square} & \CC_{0,n,\beta}\ar[d]^{h_n}\\
\YY_{0,n,\beta}^{\vec{g}}\ar[r]_{\pi_n}& \MM_{0,n,\beta}.
}
\end{equation*}

The data 
\begin{equation}\label{universal_log_twisted_curve}
(h'_n:\CC'\to \YY_{0,n,\beta}^{\vec{g}}, \{r_i=\text{order}(g_i)\}_{1\leq i\leq n},  \mathfrak{l}_n: \pi_n^*\clM_D^n\to \clM_{\clD}^n)
\end{equation}
is a log twisted curve. The universal {\em twisted} weighted curve $\CC_{0,n,\beta}^{\vec{g}}\to \YY_{0,n,\beta}^{\vec{g}}$ over $\YY_{0,n,\beta}^{\vec{g}}$ is obtained from (\ref{universal_log_twisted_curve}) by applying the construction of \cite{OLogCurv}, Section 4 (see Page \pageref{from-log-tw-to-tw-diag} for a summary).
 
 Note that over $\CC'$ there is the following diagram of log structures
\begin{eqnarray}\label{log-str-CC-incompl-diag}
\xymatrix{
h_n'^*\pi_n^*\clM_D^n\ar[r]^{h_n'^*\mathfrak{l}_n}\ar[d]_{{\pi'}_{n+1}^*h_n^b} & h_n'^* \clM_{\clD}^n\\
\pi_{n+1}^*\clM^{n+1}_D , & \\
}
\end{eqnarray}
where the vertical arrow is the  pullback of the morphism $h_n^b: h_n^*\clM_D^n\to \clM_D^{n+1}$ induced by   (\ref{preimage-of-divisor-under-univ-curve}). 
Consider the morphism 
\begin{equation}\label{h_n-stalks}
\overline{\clM}^n_{D,p}\to \overline{\clM}^{n+1}_{D,q},
\end{equation}
which is induced by   $h_n^b$ for any geometric point  $q\in \MM_{0,n+1,\beta}$ with $p=h_n(q)$. Let $C_p=h_n^{-1}(p)$. 
 The morphism (\ref{h_n-stalks}) has the following properties:

\begin{enumerate}\label{clMn-to-clMn+1-property}
\item   if $q\in C_p$ is a nodal point, i.e. $q\in D^{T_i}_{\beta_i}\cap D^{T_i\cup \{n+1\}}_{\beta_i}$ for some $i$, then, up to isomorphism,  it  is of the form $\bbN^r\to \bbN^{r+1}$, mapping $e_i$ to $e_i$, $i<r$ and $e_r$ to $e_r+e_{r+1}$;
\item  \label{clMn-to-clMn+1-property2} if $q$ is a marked point, i.e. $q\in  D_0^{\{j,n+1\}}$ for some $1\leq j\leq n$,  it  is of the form 
$\bbN^r\to \bbN^{r+1}$ mapping $e_i$ to $e_i$ for $i\leq r$; 
\item \label{clMn-to-clMn+1-property3}if $q$ is a smooth point, 
 it is of the form $\bbN^r\to \bbN^r$, mapping $e_i$ to $e_i$ for $i=1,...,r$.
\end{enumerate}

An object of $\CC_{0,n,\beta}^{\vec{g}}$ over $S$ is given by a morphism
$f: S\to\CC'$ and by a simple morphism $f^*\pi_{n+1}^*\clM^{n+1}_D \to  \clM_S'$ completing (\ref{log-str-CC-incompl-diag}) to a  commutative diagram
\begin{eqnarray}\label{univ-weigh-curv-obj-diag}
\xymatrix{
f^*h_n'^*\pi_n^*\clM_D^n\ar[r]\ar[d] & f^* h_n'^* \clM_{\clD}^n \ar[d]\\
f^*\pi_{n+1}^*\clM^{n+1}_D \ar[r] & \clM_S' \\
}
\end{eqnarray}
satisfying   conditions (\ref{from-log-tw-to-tw-properties-1}),  (\ref{from-log-tw-to-tw-properties-2}) listed on Page  \pageref{from-log-tw-to-tw-properties}.
Let $\overline{s}$ be a geometric  point of $S$ mapping to $q\in \MM_{0,n+1,\beta}$ and assume that $p=h_n(q)$ belongs to  $D^{T_1}_{\beta_1}\cap ..\cap D^{T_r}_{\beta_r}$. Consider the diagram obtained from  (\ref{univ-weigh-curv-obj-diag}) by taking the stalk at $\overline{s}$ of the associated ghost sheaves. There is a  bijection between the divisors $D^{T_i}_{\beta_i}$ containing the point  $p$ and irreducible elements of the monoid $\overline{\clM}^n_{D,\overline{s}}$.
 If  $\overline{s}$ maps to the $j$-st marked point of $C_p$, i.e.  if $q\in D_{\beta_i}^{T_1}\cap ..\cap D^{T_r}_{\beta_r}\cap D^{T_{r+1}}_{\beta_{r+1}}$, with $T_{r+1}=\{j,n+1\}$ and $\beta_{r+1}=0$,   we get the following diagram
\begin{eqnarray}\label{univ-tw-weigh-curv-obj-stalks-diag}
\xymatrix@C=2cm{
\oplus_{i=1}^r\bbN\ar[r]^-{\oplus(\times r_{T_i,\beta_i})}\ar[d] &\oplus_{i=1}^r \bbN\ar[d] \\
\oplus_{i=1}^{r+1}\bbN \ar[r]_{ \oplus(\times  \alpha_i)} & \oplus_{i=1}^{r+1}\bbN,
}
\end{eqnarray}
where the vertical arrows map $e_i$ to $e_i$, $i=1,...,r$. By  condition (\ref{from-log-tw-to-tw-properties-1}) on Page \pageref{from-log-tw-to-tw-properties}, we must have $\alpha_{i}=r_{T_i,\beta_i}$ for $i=1,..,r$ and by condition (\ref{from-log-tw-to-tw-properties-2}) on Page \pageref{from-log-tw-to-tw-properties}, $\alpha_{r+1}=r_j$, the order of the $j$-th marked point of twisted curves  parametrized by  $\YY_{0,n,\beta}^{\vec{g}}$.

Suppose that  $\overline{s}$ maps to a nodal point of $C_p$.  In this case   $q$ belongs  to  $D^{T_1}_{\beta_1}\cap ..\cap D^{T_r}_{\beta_r}\cap D^{T_{r+1}}_{\beta_{r+1}}$, where $T_{r+1}=T_r\cup\{n+1\}$ and $\beta_{r+1}=\beta_r$ for some $T_r$, $\beta_r$.   Note that we abuse the notation by using the same symbol for divisors in  $\MM_{0,n,\beta}$ and $\MM_{0,n+1,\beta}$. In this case we get  again a diagram as in (\ref{univ-tw-weigh-curv-obj-stalks-diag}),  with vertical arrows mapping  $e_i$  to $e_i$, $i<r$ and mapping $e_r$ to $e_r+e_{r+1}$. Commutativity of the diagram implies that $\alpha_r= \alpha_{r+1}=r_{T_r,\beta_r}$. 

From  the above discussion we see that    an object of $\CC_{0,n+1,\beta}^{\vec{g}}$  over $S$  encodes 
a morphism $f:S\to \MM_{0,n+1,\beta}$ and a simple morphism of log structures $f^*\clM_D^{n+1}\to \clM_S'$ as in  \cite{MO}, with coefficients $\alpha_{T,\beta'}$,  $T\subset\{1,..,n\}$, $\beta'\leq \beta$,   associated to the the irreducible components of the boundary divisor $D^T_{\beta'}$ as follows:
\begin{itemize}
\item $\alpha_{T,\beta'}=\alpha_{T\cup\{n+1\}, \beta'}=r_{T,\beta'}$,\quad $(T,\beta')\in\mathcal{I}_D$;
\item $\alpha_{\{i,n+1\},0}=r_i$,\quad $i=1,..,n$;
\end{itemize}
where the $r_{T,\beta'}$ are the integer numbers  determined as in (\ref{triple_of_integers_for_node})
for the $\beta$-admissible vector $\vec{g}$ and the class $\beta$. Observe that the collection of the  coefficients $\alpha_{T,\beta'}$ coincides with the collection of coefficients $r_{T,\beta'}$ associated to the $\beta$-admissible vector $\vec{g}\cup\{1\}$ and to the class $\beta$. In other words there is a natural morphism  $\CC_{0,n,\beta}^{\vec{g}}\to \YY_{0,n+1,\beta}^{\vec{g}\cup\{1\}}$. 

On the other hand,  it is not hard to see that  there is a morphism   $\YY_{0,n+1,\beta}^{\vec{g}\cup\{1\}}\to \CC_{0,n,\beta}^{\vec{g}}$. We start   by showing that  there is a morphism $\tilde{h}_n:\YY_{0,n+1,\beta}^{\vec{g}\cup\{1\}}\to \YY_{0,n,\beta}^{\vec{g}}$  inducing a morphism $f_{n+1}:  \YY_{0,n+1,\beta}^{\vec{g}\cup\{1\}}\to \CC'$. Consider the   universal simple morphism of log structures  $\mathfrak{l}_{n+1}: \pi_{n+1}^* \clM_D^{n+1}\to \clM_\clD^{n+1}$ over $\YY_{0,n+1,\beta}^{\vec{g}\cup\{1\}}$, where $\pi_{n+1}: \YY_{0,n+1,\beta}^{\vec{g}\cup\{1\}}\to \MM_{0,n+1, \beta}$ is the natural map, and $\clM_\clD^{n+1} $ is the log structure associated to the  divisor $\clD$ on  $\YY_{0,n+1,\beta}^{\vec{g}\cup\{1\}}$ by construction of $\YY_{0,n+1, \beta}^{\vec{g}\cup\{1\}}$. There is a natural sub-log structure $\clM_{\clD'^n}$ of $\clM_\clD^{n+1}$ which is associated to the divisor 
\begin{equation*}
\clD'^n:= \bigcup_{\stackrel{(T,\beta')\in \mathcal{I}_D}{T\subseteq[n]}} \clD_{\beta'}^{T}
\end{equation*}
The composite morphism 
$$
\pi_{n+1}^*h_n^* \clM_D^n\to \pi_{n+1}^*\clM^{n+1}_D \to \clM_\clD^{n+1}
$$
factor through $\clM_{\clD'^n}$ and $\pi_{n+1}^*h_n^* \clM_D^n\to \clM_{\clD'^n}$ is a simple morphism of log structures. This defines the morphism 
$\tilde{h}_n$.  By construction   $\clM_{\clD'^n}$   is the pullback 
of $\clM_\clD^n$ along $\tilde{h}_n$.  We conclude by observing that the  pair $(f_{n+1}, \mathfrak{l}_{n+1} )$ satisfies  diagram  (\ref{univ-weigh-curv-obj-diag}) with $f_{n+1}$ in place of $f$, giving the morphism $\YY_{0,n+1,\beta}^{\vec{g}\cup\{1\}}\to \CC^{\vec{g}}_{0,n,\beta}$. 

It is not hard to check that the two maps between $\YY_{0,n+1,\beta}^{\vec{g}\cup\{1\}}$ and $\CC_{0,n,\beta}^{\vec{g}}$ are inverse to each other.

 We use \cite{OLogCurv}, Lemma 5.3 to conclude that $\YY_{0,n+1,\beta}^{\vec{g}\cup\{1\}}$ is an open substack of $\MM_{0,n+1,\beta}^{tw}$. 
\end{pf}
\begin{rmk}
An argument similar to the proof of  Proposition \ref{family_of_twisted_curves-lem}
can be used to characterize the universal weighted twisted curve 
over $\MM_{g,n,\beta}^{tw}$ as an open substack of $\MM_{g,n+1,\beta}^{tw}$. 
\end{rmk}
Let $[f: \clC\to \clG]\in \clK_{0,n}(\clG, \beta)^{\vec{g}}$. By forgetting $f$ and keeping degrees of the restrictions of $f$ on irreducible components of $\clC$, we obtain a morphism $$\clK_{0,n}(\clG, \beta)^{\vec{g}}\to \MM_{0, n, \beta}^{tw}.$$

\begin{lem}
There exists a natural morphism $$s: \clK_{0,n}(\clG, \beta)^{\vec{g}}\to \YY_{0,n,\beta}^{\vec{g}}.$$
\end{lem}
\begin{pf}
Let $[f: (\clC, \{\sigma_i\})\to \clG]$ be an object of $\clK_{0,n}(\clG, \beta)^{\vec{g}}(\bbC)$. The pushforward $f_{*}$  defines the weight function $\clC^{gen}\to H_2^+(X, \bbZ)$. The domain $(\clC, \{\sigma_i\})$ is a genus $0$ twisted curve whose $i$-th marked gerbe ($1\leq i\leq n$) is banded by $\mu_{r_i}$. Let $x\in \clC$ be a node which separates $\clC$ into two connected components $\clC=\clC_1\cup \clC_2$. Put $f_*[\clC_1]=\beta'$ and let $T\subset [n]$ be the set of  marked points contained in $\clC_1$. Since $f$ is representable, Lemma \ref{rep_map_to_root} and equation (\ref{triple_of_integers_for_node}) imply that $x$ is a stack point of order $r_{T,\beta'}$. This defines an object of $\YY_{0,n,\beta}^{\vec{g}}$. Because indices of nodes are locally constant, extension to objects over general base schemes is straightforward.
\end{pf}


\subsection{Gerbe structures on components}\label{section:gerbe_structure}
We continue to fix a $\beta$-admissible vector $\vec{g}=(g_1,...,g_n)\in \mu_r^{\times n}$. Consider the following diagram:
\begin{eqnarray}
\label{diag_root_gerb_mod_sp}
 \xymatrix{
\clK_{0,n}(\clG,\beta)^{\vec{g}}\ar[rd]_{s}\ar@/^1pc/[rrr]^p \ar[r]_{t} & P^{\vec{g}}_n\ar[r]_{r'} \ar[d]\ar@{}[rd]|{\square} & P\ar[d]_{q'}\ar[r]_-{r}\ar@{}[rd]|{\square} & \overline{M}_{0,n}(X,\beta)\ar[d]^{q}\\
&  \mathfrak{Y}_{0,n,\beta}^{\vec{g}}\ar[r] &  \mathfrak{M}_{0,n,\beta}^{tw}\ar@{}[rd]|{\square} \ar[d]_{s'} \ar[r]& \mathfrak{M}_{0,n,\beta}\ar[d]^{s''}\\
& & \mathfrak{M}^{tw}_{0,n}\ar[r] & \mathfrak{M}_{0,n}.
}
\end{eqnarray}
Here the morphism $p: \clK_{0,n}(\clG, \beta)^{\vec{g}}\to  \overline{M}_{0,n}(X, \beta)$ is defined by sending a twisted stable map to its associated map between coarse moduli spaces. The morphism $\YY_{0,n,\beta}^{\vec{g}}\to~\MM_{0,n, \beta}^{tw}$ is defined by Proposition \ref{description_of_stack_YY}. The stacks $P$ and $P_n^{\vec{g}}$ are defined as fiber products. The morphism $t:\clK_{0,n}(\clG,\beta)^{\vec{g}}\to P^{\vec{g}}_n$ is evidently defined.


The goal of this Subsection is to prove the following

\begin{thm}\label{mod_space_to_gerbe_is_gerbe}
 Let $\clG=\sqrt[r]{\clL/X}\to X$ be an $r$-th root gerbe. Let $\vec{g}$ be a $\beta$-admissible vector for $\clG$ and a choice of $\beta\in H_2^+(X,\bbZ)$. Then the morphism $t:\clK_{0,n}(\clG,\beta)^{\vec{g}}\to P^{\vec{g}}_n$ exhibits $\clK_{0,n}(\clG,\beta)^{\vec{g}}$ as a $\mu_r$-gerbe over $ P^{\vec{g}}_n$. 
\end{thm}
\begin{pf}
We will prove that $t:\clK_{0,n}(\clG,\beta)^{\vec{g}}\to P^{\vec{g}}_n$ is a gerbe by showing that the structure morphism and the relative diagonal are epimorphisms in the sense of \cite{l-mb} Definition 3.6.  We know by Proposition \ref{unique_geom_lift} that the morphism is bijective on geometric points, hence  to  prove the first claim it is enough to show that $\clK_{0,n}(\clG,\beta)^{\vec{g}}\to P^{\vec{g}}_n$ is \'etale.  For the notion of \'etale non-representable morphisms  between algebraic stacks we refer to \cite{l-mb} Definition 4.14.
 Since the  stacks involved  are   of finite type we can use   Proposition 4.15 (ii) of \cite{l-mb}. Moreover, since we assume the stacks are noetherian, it is enough to  prove that  the lifting  criterion of \cite{l-mb}, Proposition 4.15  holds for morphisms from artinian local rings (cfr. \cite{EGAIV_IV}, 17.5.4 and \cite{EGAIV_I}, {\bf 0}, Prop. 22.1.4.)
Consider a square zero extension of artinian local rings:
\begin{equation}\label{square_0_extension}
1\to I\to B\to A\to 1.
\end{equation}
We need to show that given the outer commutative  diagram
\begin{equation*}
\xymatrix{
\Spec{A}\ar[r] \ar[d] & \clK_{0,n}(\clG,\beta)^{\vec{g}}\ar[d]\\
\Spec{B}\ar[r]\ar@{-->}[ur] & P^{\vec{g}}_n,
}
\end{equation*}
 the morphism  $\Spec{A}\to \clK_{0,n}(\clG,\beta)^{\vec{g}}$ factors 
through $\Spec{B}$.  Given a twisted stable map $\tilde{f}_A:\clC_A \to\clG$ over $\Spec{A}$ together with a lifting $\clC_A\hookrightarrow \clC_B$ of the domain curve to $\text{Spec}\, B$ and a lift $f_B:C_B\to X$ of the coarse map $f_A: C_A\to X$ to $\text{Spec}\, B$, we claim that $\tilde{f}_A$ lifts to $\Spec{B}$ uniquely. To show this first note that the exact sequence 
\begin{equation*}
H^1(I)\to H^1(\clO_{\clC_B}^*)\to H^1(\clO_{\clC_A}^*)\to H^2(I)
\end{equation*}
arising from the extension (\ref{square_0_extension}) gives an isomorphism $\Pic{\clC_A} \simeq \Pic{\clC_B}$. Indeed, $H^1(I)$ vanishes because the curves have arithmetic genus zero, and $H^2(I)$ vanishes by dimensional reasons. The morphism $\tilde{f}_A$ is defined by a line bundle $\clN_A$ satisfying $\clN_A^{\otimes r}\simeq \pi_A^*f_A^*\mathcal{L}$, where $\pi_A:\clC_A\to C_A$ is the map to the coarse curve. The isomorphism $\Pic{\clC_A}\simeq \Pic{\clC_B}$ yields a line bundle $\clN_B$  on $\clC_B$ which satisfies $\clN_B^r\simeq \pi_B^*f_B^*\mathcal{L}$. (Again $\pi_B: \clC_B\to C_B$ denotes the map to the coarse curve.)  This defines the desired extension $\tilde{f}_B:\clC_B\to \clG$. 

It remains to show that the relative diagonal morphism $$\clK_{0,n}(\clG,\beta)^{\vec{g}}\to \clK_{0,n}(\clG,\beta)^{\vec{g}}\times_{P^{\vec{g}}_n}\clK_{0,n}(\clG,\beta)^{\vec{g}}$$ is locally surjective, namely that any two local sections are locally isomorphic. Let $T\to P_n^{\vec{g}}$ be a morphism which gives rise the twisted curve $\clC_T$ over $T$ with coarse curve $\pi_T: \clC_T\to C_T$ and stable map $f_T: C_T\to X$. Consider the base-change 
\begin{equation}\label{pullback_to_T} 
\clK_{0,n}(\clG,\beta)^{\vec{g}}\times_{P_n^{\vec{g}}} T\to T.
\end{equation}
Sections of (\ref{pullback_to_T}) are twisted stable maps $\clC_T\to \clG$ which induce $f_T$. Such maps are defined by line bundles $\clN$ on $\clC_T$ satisfying $\clN^{\otimes r}\simeq \pi_T^* f_T^* \mathcal{L}$. Two sections are isomorphic if and only if their defining line bundles are isomorphic. 

Let  $f$ and $f'$ be two sections of (\ref{pullback_to_T}) with line bundles $\clN$ and $\clN'$. Note that the line bundle $\clN\otimes \clN'^{\vee}$ 
is a pullback from the coarse moduli space $C_T$ since it carries 
trivial representations of the special points of any geometric fiber of  $\clC_T$. This is due to the fact that   $\clN$ and  $\clN'$ induce twisted stable maps 
with same $\beta$-admissible vector. Note moreover that   $(\clN\otimes \clN'^{\vee})^{\otimes r}$is trivial:
$$(\clN\otimes \clN'^{\vee})^{\otimes r}\simeq \pi_T^*f_T^*\mathcal{L}\otimes \pi_T^*f_T^*\mathcal{L}^{\vee}\simeq \clO_{\clC_T}.$$
As a consequence,  $\clN\otimes \clN'^{\vee}$
restricts to a trivial  line bundle over any irreducible component of any geometric fiber of $\clC_{T}$.
By the Theorem on Cohomology and Base change, $\clN\otimes \clN'^{\vee}$ is isomorphic to the pullback of a line bundle from the base $T$. Every line  bundle  over  a  scheme is locally trivial, therefore up to base change by an \'etale morphism $T'\to T$ the two line bundles $\clN$ and $\clN'$ are isomorphic. 

For any object $\xi$ of $\clK_{0,n}(\clG,\beta)^{\vec{g}}$ over some $x\in P^{\vec{g}}_n(T)$ the sheaf of relative automorphisms  $Aut_x(\xi)$ is a sheaf of abelian groups. Therefore $\clK_{0,n}(\clG,\beta)^{\vec{g}}\to P^{\vec{g}}_n$ is a gerbe banded by a sheaf of abelian groups (\cite{Gir}, Proposition IV 1.2.3 (i)). For any $\xi$ in $\clK_{0,n}(\clG,\beta)^{\vec{g}}$,  there is a natural identification $Aut_x(\xi)\simeq (\mu_r)_T$. Indeed automorphisms of $\xi$ leaving $x$ fixed are automorphisms of a line bundle over a twisted curve $p: \clC\to T$ whose $r$-th power is the identity. The claim follows since for a family of twisted curves $p_*\clO_\clC\simeq \clO_T$. Such a collection of natural identifications is compatible with restrictions and isomorphisms. This means by definition that the $\clK_{0,n}(\clG,\beta)^{\vec{g}}$ is banded by  $(\mu_r)_{P^{\vec{g}}_n}$ (\cite{Gir}, IV, Definition 2.2.2).
\end{pf}



\subsection{Root gerbe structures on components}\label{section:root_gerbe_structure}
Consider the following diagram
\begin{equation}\label{big_square1}
\xymatrix{
P_{n+1}^{\vec{g}\cup \{1\}}\ar[d]^{v}\ar[r] \ar@{}[rd]|{\square}& \overline{M}_{0,n+1}(X, \beta)\ar[d]\ar[r]\ar@{}[rd]|{\square} & \overline{M}_{0,n}(X,\beta)\ar[d] \\
\YY_{0,n+1, \beta}^{\vec{g}\cup\{1\}}\ar[r] & \MM_{0,n+1,\beta} \ar[r] & \MM_{0,n,\beta}. 
}
\end{equation}
The stack $P_{n+1}^{\vec{g}\cup\{1\}}$ is defined by the cartesian square on the left. The square on the right is cartesian by Proposition \ref{universal_family}. The existence of the morphism $\YY_{0,n+1,\beta}^{\vec{g}\cup \{1\}}\to \YY_{0,n,\beta}^{\vec{g}}$ implies that the outer square in (\ref{big_square1}) is equivalent to the outer square of the following diagram
\begin{equation}\label{big_square2}
\xymatrix{
P_{n+1}^{\vec{g}\cup \{1\}}\ar[d]^{v}\ar[r]^{\phi}& P_n^{\vec{g}}\ar[d]\ar[r]\ar@{}[rd]|{\square}& \overline{M}_{0,n}(X, \beta)\ar[d]\\
\YY_{0,n+1,\beta}^{\vec{g}\cup\{1\}}\ar[r] &\YY_{0,n,\beta}^{\vec{g}}\ar[r]  &\MM_{0,n,\beta}.
}
\end{equation}
Since (\ref{big_square1}) and the right side of (\ref{big_square2}) are cartesian, the left part of (\ref{big_square2}) is also cartesian.

By construction the stack $P_n^{\vec{g}}$ parametrizes the data $$((\clC, \sigma_1,...,\sigma_n), (f: (C, p_1,...,p_n)\to X))$$ where $(\clC, \sigma_1,...,\sigma_n)$ is an $n$-pointed twisted curve with isotropy group at $\sigma_i$ being $\mu_{r_i}$, $(C, p_1,...,p_n)$ is the coarse curve, and $[f]\in \overline{M}_{0,n}(X, \beta)$ is a stable map.  The universal twisted curve over $P_n^{\vec{g}}$ is given by the morphism $P_{n+1}^{\vec{g}\cup \{1\}}\to  P_n^{\vec{g}}$, and the universal stable map is obtained from the composition $$u: P_{n+1}^{\vec{g}\cup \{1\}}\to \overline{M}_{0,n+1}(X, \beta)\overset{ev_{n+1}}{\longrightarrow} X.$$ 

The purpose of this Subsection is to prove the following refinement of Theorem \ref{mod_space_to_gerbe_is_gerbe}:

\begin{thm}\label{mod_space_to_root_gerbe_is_root_gerbe}
 $\clK_{0,n}(\clG,\beta)^{\vec{g}}$ is a root  gerbe over the stack $P^{\vec{g}}_n$.
\end{thm}
\begin{pf}
We construct a line bundle over $P_n^{\vec{g}}$ and show that the stack of its $r$-th roots admits a representable morphism to $\clK_{0,n}(\clG, \beta)^{\vec{g}}$ which covers the identity map on $P_n^{\vec{g}}$.

Recall that for a $\beta$-admissible vector $\vec{g}=(g_1,...,g_n)$ we have defined triples $(\rho_i, r_i, m_i), 1\leq i\leq n$ of integers in (\ref{the_triple_of_integers}). For $1\leq i\leq n$ we choose $d_i\in \bbZ$ such that 
\begin{equation}\label{requirement_for_d_i}
g_i=\exp(\frac{2\pi \sqrt{-1}}{r_i}d_i), \quad \text{and } \sum_{i=1}^n \frac{d_i}{r_i}=\frac{1}{r}\int_{\beta}c_1(\clL).
\end{equation}
This is possible because of (\ref{phase_shift_eq}). Note that $d_i$ depends on $\beta$.  

Associated to a pair $(T, \beta')$ which indexes a boundary divisor of $\MM_{0,n,\beta}$, we have defined a triple $(\theta_{T, \beta'}, r_{T, \beta'}, m_{T, \beta'})$ of integers in (\ref{triple_of_integers_for_node}). We define another integer $d_{T, \beta'}$ such that 
\begin{equation}\label{requirement_for_d_node}
\sum_{i\in T} \frac{d_i}{r_i} +\frac{d_{T, \beta'}}{r_{T, \beta'}}=\frac{1}{r}\int_{\beta'} c_1(\clL).
\end{equation}
Note that (\ref{requirement_for_d_node}) implies $d_{T, \beta'}=-d_{T^C, \beta-\beta'}$ and $\langle\frac{d_{T, \beta'}}{r_{T, \beta'}} \rangle=\theta_{T, \beta'}$. 

For $1\leq i\leq n$, let $S_i\subset \YY_{0,n+1, \beta}^{\vec{g}\cup \{1\}}$ denote the pullback of the  $i$-th marked section divisor from $\MM_{0,n+1,\beta}$. Define a line bundle over $\YY_{0,n+1, \beta}^{\vec{g}\cup \{1\}}$ as follows:
\begin{equation}\label{L_YY_eq}
\clL_\YY:=\mathcal{O}_{\YY_{0,n+1, \beta}^{\vec{g}\cup \{1\}}}\left(\sum_{1\leq i\leq n} \frac{d_i}{r_i} S_i -\sum_{(T, \beta')\in \mathcal{I}_D} \frac{d_{T, \beta'}}{r_{T, \beta'}}D_{\beta'}^{T\cup\{n+1\}}\right).
\end{equation}
Here $\mathcal{I}_D$ is the set of pairs $(T, \beta')$ defined on page \pageref{calI_D_def_intext}.

\begin{lem}\label{defining_L_2_lem}
There exists a line bundle $\clL_2$ over $P_n^{\vec{g}}$ such that 
\begin{equation}\label{defining_L_2}
(v^*\clL_\YY)^{\otimes r}\otimes (u^* \clL)^{-1}\simeq \phi^* \clL_2.
\end{equation}
\end{lem}
\begin{pf}
As in \cite{BC}, it suffices to check that the degree of the line bundle $(v^*\clL_\YY)^{\otimes r}\otimes (u^* \clL)^{-1}$ restricted to any component of any fiber of $\phi$ is zero. The argument works because $(v^*\clL_\YY)^{\otimes r}\otimes (u^* \clL)^{-1}$ is in fact a pullback from $\MM_{0,n+1,\beta}$.
Indeed $\clL$ is a line bundle over a scheme and $(v^*\clL_\YY)^{\otimes r}$
carries trivial representations of the automorphisms groups of stacky points of $\YY_{0,n+1, \beta}^{\vec{g}\cup \{1\}}$  relative to $\MM_{0,n+1,\beta}$. 
Let $\clC$ be a fiber of $\phi$, with coarse curve $C$. Denote by $f: C\to X$ the corresponding stable map to $X$. Let $\clC^0\subset \clC$ be an irreducible component with coarse curve $C^0$. Let $x_1, ..., x_m$ be nodes of $\clC$ that are contained in $\clC^0$. Let $T_j\subset [n], 1\leq j\leq m$ be the marked points contained in the subcurves $\clC^j\subset \clC$ which are connected to $\clC^0$ at $x_j$, and $T_0$ the marked points contained in $\clC^0$. Then $[n]=T_0\cup T_1\cup ...\cup T_m$. Put $\beta_0:=f_*[C^0]$ and $\beta_j:=f_*[C^j]$ (here $C^j$ is the coarse curve of $\clC^j$). 

We need some properties about restrictions of these line bundles.\\
{\bf Claim.} Consider the line bundle $\clL^{(T,\beta')}:=\clO_{\YY_{0,n+1,\beta}^{\vec{g}\cup\{1\}}}(\frac{1}{r_{T,\beta'}}D_{\beta'}^{T\cup\{n+1\}})$. Let  $\clC$ be  a geometric fiber of $\YY_{0,n+1,\beta}^{\vec{g}\cup\{1\}}\to \YY_{0,n,\beta}^{\vec{g}}$. 
\begin{itemize}
\item If there is no node $e$ in $\clC$ such that the two connected components of the normalization of $\clC$ at $e$  have degrees $\beta'$, resp. $\beta''$ (such that $\beta'+\beta''=\beta$),  and contain marked points with indices in $T$, resp. $T^C$,  then $\clL^{(T,\beta')}|_\clC$ is trivial.
\item If there is such a node, let $\clC_1$ and $\clC_2$ be  the two connected components of the partial normalization at the node. Suppose that the preimages of the marked gerbes with indices in $T$ are contained in  $\clC_1$.
  Let $\overline{\clC}_1\subset \clC_1$ and $\overline{\clC}_2\subset \clC_2$ be the two irreducible components in $\clC_1$ and $\clC_2$ containing the node $e$. Then
\begin{enumerate}
\item \label{clLC1} $\clL^{(T,\beta')}|_{\overline{\clC}_1}\simeq \clO_{\overline{\clC}_1}(-\frac{1}{r_{T,\beta'}})$;
\item   \label{clLC2}  $\clL^{(T,\beta')}|_{\overline{\clC}_2}\simeq \clO_{\overline{\clC}_2}(\frac{1}{r_{T,\beta'}})$. 
\item \label{clLC'} the restriction of $\clL^{(T, \beta')}$ to any other component $\clC'$ of $\clC$ is trivial. 
\end{enumerate}
\end{itemize}
{\em Proof of Claim.}
The first property follows since $\clC$ misses $D_{\beta'}^{T\cup \{n+1\}}$. For the second property, first note that $\clC_1=\clC\cap D^{T\cup\{n+1\}}_{\beta'}$ and $\clC_2=\clC\cap D^{T^C}_{\beta''}$.  Therefore  $\clC_2\cap D^{T\cup\{n+1\}}_{\beta'}$ is one point and (\ref{clLC2}) above follows.  (\ref{clLC'}) for any component $\clC'\subset \clC_2$ also follows. Moreover, since $\mathcal{O}(\frac{1}{r_{T,\beta'}}(D^{T}_{\beta'}\cup  D^{T\cup\{n+1\}}_{\beta'}))$ over $\YY_{0, n+1, \beta}^{\vec{g}\cup\{1\}}$ is the pullback of  $\mathcal{O}(\frac{1}{r_{T,\beta'}}D^{T}_{\beta'})$ over $\YY_{0, n, \beta}^{\vec{g}}$,  its restriction to geometric fibers has to be trivial. This implies that $\clO_{\YY_{0, n+1, \beta}^{\vec{g}\cup \{1\}}}(\frac{1}{r_{T,\beta'}} D^T_{\beta'})|_{\overline{\clC}_2}\simeq \clO_{\overline{\clC}_2}(-\frac{1}{r_{T,\beta'}} )$. 
Statements (\ref{clLC1}) and (\ref{clLC'}) for $\clC'\subset \clC_1$ follow by symmetry.$\ \spadesuit$ 

Notice that the following diagram
\begin{eqnarray}
\label{P_frY_cart_diag}
\xymatrix{
P_{n+1}^{\vec{g}\cup\{1\}}\ar[d]\ar[r]^v & \YY_{0,n+1, \beta}^{\vec{g}\cup\{1\}}\ar[d]\\
P^{\vec{g}}_n\ar[r] & \YY^{\vec{g}}_{0, n, \beta}\\
}
\end{eqnarray}
is cartesian, hence the fibers of $P_{n+1}^{\vec{g}\cup\{1\}}\to P^{\vec{g}}_n$ are isomorphic to the fibers of $ \YY_{0,n+1, \beta}^{\vec{g}\cup\{1\}}\to \YY^{\vec{g}}_{0, n, \beta}$. 
Applying the above  Claim and the fact that $d_{T, \beta'}=-d_{T^C, \beta-\beta'}$, we find
 \begin{eqnarray} \label{L_YY_restr_eq}
v^*\clL_\YY|_{\clC^0}=\mathcal{O}_{\clC^0}(\sum_{i\in T_0} \frac{d_i}{r_i}S_i-\sum_{1\leq j\leq m}\frac{d_{T_j, \beta_j}}{r_{T_j, \beta_j}}x_j).
\end{eqnarray}
By (\ref{requirement_for_d_node}) we have $$\sum_{i\in T_j}\frac{d_i}{r_i}+\frac{d_{T_j, \beta_j}}{r_{T_j, \beta_j}}=\frac{1}{r}\int_{\beta_j}c_1(\clL).$$
By (\ref{requirement_for_d_i}) we have $$\sum_{i\in T_0}\frac{d_i}{r_i}+\sum_{1\leq j\leq m}\sum_{i\in T_j}\frac{d_i}{r_i} =\frac{1}{r}\int_{\beta} c_1(\clL).$$
Since $\beta=\beta_0+\sum_{1\leq j\leq m}\beta_j$, we find that the degree of $v^*\clL_\YY|_{\clC^0}$ is 
\begin{equation*}
\begin{split}
&\sum_{i\in T_0} \frac{d_i}{r_i}-\sum_{1\leq j\leq m} \frac{d_{T_j, \beta_j}}{r_{T_j, \beta_j}}\\
=&(\frac{1}{r}\int_\beta c_1(\clL)- \sum_{1\leq j\leq m}\sum_{i\in T_j}\frac{d_i}{r_i}) -\sum_{1\leq j\leq m}(\frac{1}{r}\int_{\beta_j}c_1(\clL)-\sum_{i\in T_j}\frac{d_i}{r_i})\\
=&\frac{1}{r}\int_{\beta_0} c_1(\clL).
\end{split}
\end{equation*}
Thus the degree of $(v^*\clL_\YY)^{\otimes r}\otimes (u^* \clL)^{-1}|_{\clC^0}$ is zero, as desired. 
\end{pf}

Let $\mathcal{P}_n^{\vec{g}}=\sqrt[r]{\clL_2/P_n^{\vec{g}}}$ be the stack of $r$-th roots of $\clL_2$. Denote by $\pi_{P_n}: \mathcal{P}_n^{\vec{g}} \to P_n^{\vec{g}}$. Consider the pull-back of $P_{n+1}^{\vec{g}\cup\{1\}}\to P_n^{\vec{g}}$ to $\mathcal{P}_n^{\vec{g}}$ via $\pi_{P_n}$:
\begin{equation}\label{pulling_back_family_of_twisted_curves}
\xymatrix{
\mathcal{P}_{n+1}^{\vec{g}\cup \{1\}}\ar[d]^{\phi'}\ar[r]^{\pi_{P_{n+1}}}& P_{n+1}^{\vec{g}\cup \{1\}}\ar[d]^{\phi}\ar[r]^u & X\\
 \mathcal{P}_n^{\vec{g}}\ar[r]^{\pi_{P_n}} & P_n^{\vec{g}}.
}
\end{equation}
\begin{lem}\label{lem:induced_family_of_twisted_stable_map}
There exists a family of twisted stable maps 
\begin{equation}\label{induced_family_of_twisted_stable_map}
\xymatrix{
\mathcal{P}_{n+1}^{\vec{g}\cup \{1\}}\ar[d] \ar[r] & \sqrt[r]{\clL/X}\\
\mathcal{P}_n^{\vec{g}}. 
}
\end{equation}
\end{lem}
\begin{pf}
Let $\clL_2^{1/r}$ be the universal $r$-th root line bundle of $\pi_{P_n}^*\clL_2$ over $\mathcal{P}_n^{\vec{g}}$. In (\ref{pulling_back_family_of_twisted_curves}) we calculate
\begin{equation*}
\begin{split}
\pi_{P_{n+1}}^*u^*\clL \simeq &\pi_{P_{n+1}}^*((v^*\clL_\YY)^{\otimes r}\otimes (\phi^* \clL_2)^{-1}) \quad \text{by } (\ref{defining_L_2})\\
\simeq & (\pi_{P_{n+1}}^*v^*\clL_\YY)^{\otimes r}\otimes (\phi'^*\pi_{P_{n}}^*\clL_2)^{-1}\\
\simeq & (\pi_{P_{n+1}}^*v^*\clL_\YY \otimes \phi'^*(\clL_2^{1/r})^{-1})^{\otimes r}.
\end{split}
\end{equation*}
The line bundle $\clL_\clP:=\pi_{P_{n+1}}^*v^*\clL_\YY \otimes \phi'^*(\clL_2^{1/r})^{-1}$ defines a morphism 
\begin{equation}\label{map_given_by_line_bundle}
\mathcal{P}_{n+1}^{\vec{g}\cup\{1\}}\to \sqrt[r]{\clL/X}.
\end{equation}
Since $\phi': \mathcal{P}_{n+1}^{\vec{g}\cup\{1\}}\to \mathcal{P}_n^{\vec{g}}$ is a family of twisted curves, to show that the morphism (\ref{map_given_by_line_bundle}) is a family of twisted stable maps we need to check that it is representable. For this purpose it suffices to work on geometric fibers of $\phi'$. 

Let $(\clC, \{\sigma_i\})$ be a geometric fiber of $\phi'$, with coarse curve $(C, \{p_i\})$ and a given stable map $\bar{f}: (C, \{p_i\})\to X$ of degree $\beta$. The restriction of (\ref{map_given_by_line_bundle}) to $\clC$ is given by the line bundle $\clL_\clP|_{\clC}$, and fits into the following diagram:
\begin{equation*}
\xymatrix{
\clC\ar[r]^-{\tilde{f}}\ar[d] & \sqrt[r]{\clL/X} \ar[d]\\
C\ar[r]^-{\bar{f}} &X.
}
\end{equation*}
By the choices of $d_i, 1\leq i\leq n$ as in (\ref{requirement_for_d_i}), the action of the stabilizer group of $\sigma_i$ on $\clL_\clP|_{\sigma_i}$ is given by the element $g_i\in \mu_r$. This is due to the fact that by construction
$Aut(\sigma_i)$ only  acts non trivially  on the fibers of $\clL_\YY$.
By the choice of $r_{T, \beta'}, d_{T, \beta'}$ as in (\ref{triple_of_integers_for_node}) and (\ref{requirement_for_d_node}), on the restriction of $\clL_\clP|_{\clC}$ to a node $x\in\clC$, the action of the stabilizer group of $x$ is given by $\exp(2\pi\sqrt{-1}\frac{d_{T,\beta'}}{r_{T, \beta'}})=\exp(2\pi\sqrt{-1}\frac{m_{T, \beta'}}{r_{T, \beta'}})$ (when viewed $\clL_\clP|_{\clC}$ as a line bundle on one of the branches meeting at $x$). Since by construction $m_{T, \beta'}$ and $r_{T, \beta'}$ are co-prime, it follows from Lemma \ref{rep_map_to_root} that the restriction of $\tilde{f}$ to any irreducible component of $\clC$ is representable. Therefore $\tilde{f}$ is representable. Hence (\ref{map_given_by_line_bundle}) gives the family (\ref{induced_family_of_twisted_stable_map}) of twisted stable maps we want.
\end{pf}

As observed in the proof of Lemma \ref{lem:induced_family_of_twisted_stable_map}, the family (\ref{induced_family_of_twisted_stable_map}) induces a morphism $\mathcal{P}_{n}^{\vec{g}}\to \clK_{0,n}(\clG, \beta)$. The choices of $d_i$ ensures that actions of isotropy groups at marked points are given by the $\beta$-admissible vector $\vec{g}$. In other words (\ref{induced_family_of_twisted_stable_map}) gives a morphism 
\begin{equation}\label{map_from_gerbe_to_moduli}
\mathcal{P}_{n}^{\vec{g}}\to \clK_{0,n}(\clG, \beta)^{\vec{g}}, 
\end{equation}
 which fits into the following diagram:
\begin{equation}\label{gerbe_and_moduli_space}
\xymatrix{
\mathcal{P}_{n}^{\vec{g}}\ar[r]\ar[d] & \clK_{0,n}(\clG, \beta)^{\vec{g}}\ar[d]\\
P_n^{\vec{g}}\ar[r]^{id} & P_n^{\vec{g}}.
}
\end{equation}
It is straightforward to check that the diagram is commutative. Indeed, an object
$g:T\to \mathcal{P}_n^{\vec{g}}$  is given by  
\begin{equation}\label{object_of_clPn}
\pi_g:\clC\to T,\quad  f:C\to X,\quad  \clN,
\end{equation}
 where  $\pi_g:\clC\to T$ is a  family of twisted curves, $f:C\to X$ is a stable map from the coarse moduli space of  $\clC$ to $X$, and $\clN$ is a line bundle over $T$ together with an isomorphism $\clN^{\otimes r}\simeq g^*\pi_{P_n}^*\clL_2$. Note that $\clC\to T$ is obtained as the pull-back of $\phi': \mathcal{P}_{0,n+1, \beta}^{\vec{g}\cup\{1\}}\to \mathcal{P}_{0,n,\beta}^{\vec{g}}$ via $g: T\to \mathcal{P}_{0,n,\beta}^{\vec{g}}$, hence there is a morphism $g_\clC:\clC\to \mathcal{P}_{0,n+1,\beta}^{\vec{g}\cup\{1\}}$.
 
 In (\ref{gerbe_and_moduli_space}), the morphism  $\mathcal{P}_{n}^{\vec{g}}\to P_n^{\vec{g}}$ just forgets the line bundle $\clN$. The morphism $\mathcal{P}_{n}^{\vec{g}}\to \clK_{0,n}(\clG, \beta)^{\vec{g}}$ takes an object (\ref{object_of_clPn}) to a twisted stable map $\clC\to \clG$ defined by the line bundle  $g_\clC^*v^*\clL_\YY \otimes \pi_g^*\clN^{-1}$. The morphism $\clK_{0,n}(\clG, \beta)^{\vec{g}}\to P^{\vec{g}}_n$ retains $\pi_g: \clC\to T$ and $f:C\to X$. Hence (\ref{gerbe_and_moduli_space}) is commutative.
 
An easy analysis on automorphism groups of (\ref{object_of_clPn}) and twisted stable maps shows that the morphism (\ref{map_from_gerbe_to_moduli}) is representable. By \cite{Gir}, IV, Proposition 2.2.6, a representable morphism between two gerbes banded by the same group is an isomorphism. Since both $\mathcal{P}^{\vec{g}}_n$ and $\clK_{0,n}(\clG,\beta)^{\vec{g}}$ are gerbes banded by the group $\mu_r$ (c.f. Theorem \ref{mod_space_to_gerbe_is_gerbe}), (\ref{map_from_gerbe_to_moduli}) is an isomorphism. This completes the proof of Theorem \ref{mod_space_to_root_gerbe_is_root_gerbe}.
\end{pf}

\section{Gromov-Witten invariants}\label{GW_theory}

\subsection{Virtual Fundamental Class}
The moduli space $\overline{M}_{0,n}(X,\beta)$ has a perfect obstruction theory  relative to $\MM_{0,n}$  which is given by  
\begin{equation*}
E^{\bullet}:=R\pi_*(f^*\Omega_X\otimes \omega_{\pi})\to  L_{ \overline{M}_{0,n}(X,\beta)/\mathfrak{M}_{0,n}},
\end{equation*}
where 
\begin{equation*}
\xymatrix{
C\ar[d]^{\pi}\ar[r]^{f} & X\\
\overline{M}_{0,n}(X,\beta)
}
\end{equation*}
is the universal stable map. Notice that since the morphism $s''$ in diagram (\ref{diag_root_gerb_mod_sp}) is \'etale, $E^{\bullet}$ is also a perfect obstruction theory relative to $\mathfrak{M}_{0,n,\beta}$.

Consider the universal twisted stable map to the gerbe $\clG$:
\begin{equation*}
\xymatrix{
\clC\ar[d]^{\tilde{\pi}}\ar[r]^{\tilde{f}}& \clG\\
\clK_{0,n}(\clG,\beta).
}
\end{equation*}
 According to \cite{AGV2} the moduli stack
$\clK_{0,n}(\clG,\beta)$ has a perfect obstruction theory $\tilde{E}^\bullet$ relative to $\MM^{tw}_{0,n}$ given by
\begin{equation*}
\tilde{E}^{\bullet}:=R\tilde{\pi}_*(\tilde{f}^*\Omega_{\clG}\otimes \omega_{\tilde{\pi}})\to  L_{ \clK_{0,n}(\clG,\beta)/\mathfrak{M}^{tw}_{0,n}}.
\end{equation*}
Since the morphism $\YY_{0,n,\beta}^{\vec{g}}\to \MM_{0,n,\beta}^{tw}$ is \'etale (Lemma \ref{some_map_is_etale}), we can view $\tilde{E}^\bullet$ as a perfect obstruction theory relative to $\YY_{0,n,\beta}^{\vec{g}}$. 

The complex $\tilde{E}^{\bullet}$  turns out to be the pullback of $E^{\bullet}$ as an object in $\clD_{Coh}( \clK_{0,n}(\clG,\beta))$.


\begin{lem}\label{PROTHtwvsnontw}
 There is a natural isomorphism of objects in $\clD_{Coh}( \clK_{0,n}(\clG,\beta))$,
\begin{equation*}
p^*E^\bullet \isomto \tilde{E}^\bullet.
\end{equation*}
\end{lem}

\begin{pf}
We will prove the statement for  $\tilde{E}^{\vee\bullet}=R\tilde{\pi}_*\tilde{f}^*T_\clG$ and $E^{\vee\bullet}=R\pi_*f^*T_X$. Consider the complex $Lp^*R\pi_*f^*T_X$ in $\clD_{coh}(\clK_{0,n}(\clG,\beta))$. It suffices to show that $Lp^*R\pi_*f^*T_X\simeq R\tilde{\pi}_*\tilde{f}^*T_\clG$. For this we consider the diagram 
\begin{equation*}
\xymatrix{
\clC\ar[dr]^\rho\ar[ddr]^{\tilde{\pi}}\ar[drr]\ar[rr]^{\tilde{f}} & & \clG\ar[dr]^\epsilon\\
&p^*C\ar[d]^{\pi_1}\ar[r]^{p_1} &C\ar[d]^\pi\ar[r]^f & X\\
&\clK_{g,n}(\clG, \beta)\ar[r]^p & \overline{M}_{0,n}(X, \beta).
}
\end{equation*}
Observe that $\epsilon^*T_X\simeq T_\clG$. Also we have $R\rho_*L\rho^*\simeq Id$ because the map $\rho$ is the relative coarse moduli space for the map $\tilde{\pi}$. The  arrow $p$ is flat, as follows from Remark  \ref{YY_to_MM-flat} and from the gerbe structure of $\clK_{g,n}(\clG, \beta)^{\vec{g}}$ over $P_n^{\vec{g}}$ for any $\beta$-admissible vector $\vec{g}$. The arrow $\pi$ is flat because it  is the structure morphism of the universal curve. 
 Moreover  the square in the  above diagram is cartesian, hence  
 we calculate (based on  \cite{l-mb}, Proposition 13.1.9)
\begin{equation*}
\begin{split}
Lp^*R\pi_*f^*T_X& \simeq R\pi_{1*}Lp_1^*f^*T_X\\
&\simeq R\pi_{1*}R\rho_*L\rho^*Lp_1^*f^*T_X\\
&\simeq R\pi_{1*}R\rho_*\tilde{f}^*\epsilon^*T_X\\
&\simeq R\pi_{1*}R\rho_*\tilde{f}^*T_\clG\\
&\simeq R\tilde{\pi}_*\tilde{f}^*T_\clG.
\end{split}
\end{equation*}
Since $p^*$  is  exact,  we write $p^*$ for $Lp^*$. 
\end{pf}

\begin{rmk}
The composite  morphism $p^*E^\bullet \to p^*L_{\overline{M}_{0,n}(X,\beta)/\MM_{0,n}}\stackrel{\sim}{\to} L_{\clK_{0,n}(\clG,\beta)^{\vec{g}}/\MM_{0,n}^{tw}}$
in $D_{Coh}(\clK_{0,n}(\clG,\beta)^{\vec{g}})$  is the same as  $E^{\vee\bullet} \to  L_{\clK_{0,n}(\clG,\beta)^{\vec{g}}/\MM_{0,n}^{tw}}$. This follows from functorial properties  of the  cotangent complex (\cite{Ill_I} Ch.2  Sect. 1 and 2). The relative obstruction theories  $E^\bullet$ and  $E^{\vee\bullet}$ are built from functorial morphisms between the cotangent complexes of the target scheme or stack  and of the universal objects (cfr. \cite{BF} Sect. 6).
\end{rmk}

\begin{thm}\label{pushforward_vir}
Let $\clG\to X$ an $r$-th root of a line bundle. Let $\vec{g}$ be a $\beta$-admissible vector. Then
the following relation between the virtual fundamental classes holds: 
\begin{eqnarray}
\label{virt_fund_class_eq}
p_*[\clK_{0,n}(\clG,\beta)^{\vec{g}}]^{vir}=\frac{1}{r}[ \overline{M}_{0,n}(X,\beta)]^{vir}.
\end{eqnarray}
\end{thm}
\begin{pf}
Consider the diagram (\ref{diag_root_gerb_mod_sp}) again. Observe that the following hold:
\begin{itemize}
\item $\mathfrak{M}_{0,n,\beta}$ and  $\mathfrak{Y}_{0,n,\beta}^{\vec{g}}$ are smooth Artin stacks of the same pure dimension;
\item  The morphism  $\mathfrak{Y}_{0,n,\beta}^{\vec{g}}\to \mathfrak{M}_{0,n,\beta}$ is of Deligne-Mumford type and of pure degree;
\item The morphism $r'\circ r$ is proper (because being proper is preserved by base change);
\item  $\overline{M}_{0,n}(X,\beta)\to  \mathfrak{M}_{0,n,\beta}$ has a perfect relative obstruction theory $E^{\bullet}$ inducing a perfect relative obstruction theory on
 $P_n^{\vec{g}}\to  \mathfrak{Y}_{0,n,\beta}^{\vec{g}}$. 
\end{itemize} 
Therefore we can apply Theorem 5.0.1 of \cite{Cost03} and conclude that $(r'\circ r)_*[P_n^{\vec{g}}]^{vir}=[ \overline{M}_{0,n}(X,\beta)]^{vir}$, where the multiplicative factor is 1 because by construction (\cite{MO}, Theorem 4.1) $\mathfrak{Y}_{0,n,\beta}^{\vec{g}}\to \mathfrak{M}_{0,n,\beta}$ is an isomorphism outside a locus of codimension $1$,  hence its base change has virtual degree equal to 1. Since $t$ is \'etale $[\clK_{0,n}(\clG,\beta)^{\vec{g}}]^{vir}=deg(t)\cdot [P_n^{\vec{g}}]^{vir}$. Since, by Theorem \ref{mod_space_to_gerbe_is_gerbe} or \ref{mod_space_to_root_gerbe_is_root_gerbe}, $t$ is of degree $\frac{1}{r}$, (\ref{virt_fund_class_eq}) follows. 
\end{pf}


\subsection{Genus $0$ invariants}
Let
$$\epsilon: \clG=\sqrt[r]{\mathcal{L}/X} \to X$$ 
be a $\mu_{r}$-root gerbe.
Then the inertia stack admits the following decomposition
$$I\clG=\coprod_{g\in \mu_{r}}\clG_{g},$$
where $\clG_{g}$ is a  root gerbe isomorphic to  $\clG$. Let $\epsilon_{g}:\clG_{g}\to X$
be the induced morphism.  On each component there is an isomorphism between the rational cohomology groups 
$$\epsilon_{g}^{*}: H^{*}(X,\mathbb{Q})\overset{\simeq}{\longrightarrow} H^{*}(\clG_{g}, \mathbb{Q}).$$

Let $\vec{g}=(g_1,...,g_n)$ be a $\beta$-admissible vector. There are evaluation maps $$ev_i: \mathcal{K}_{0,n}(\sqrt[r]{\mathcal{L}/X},\beta)^{\vec{g}}\to \bar{I}(\clG)_{g_i},$$
where $\bar{I}(\clG)_{g_i}$ is a component of  the  {\em rigidified inertia stack} $\bar{I}(\clG)=\cup_{g\in \mu_r} \bar{I}(\clG)_g$. Although the evaluation maps  $ev_i$ do not take values in  $I\clG$, as explained in \cite{AGV2}, Section 6.1.3, one can still define a pull-back map at cohomology level,
$$ev_i^*: H^*(\clG_{g_i}, \mathbb{Q})\to H^*(\mathcal{K}_{0,n}(\sqrt[r]{\mathcal{L}/X},\beta)^{\vec{g}}, \mathbb{Q}).$$
Given $\delta_i\in H^*(\clG_{g_{i}},\mathbb{Q})$ for $1\leq i\leq n$ and   integers $k_i\geq 0,1\leq i\leq n$,  one can define descendant orbifold Gromov-Witten invariants 
$$\langle \delta_1\bar{\psi}_1^{k_1},\cdots, \delta_n\bar{\psi}_n^{k_n}\rangle_{0, n, \beta}^\clG:=\int_{[\mathcal{K}_{0,n}(\sqrt[r]{\mathcal{L}/X},\beta)^{\vec{g}}]^{vir}}\prod_{i=1}^n ev_i^*(\delta_i)\bar{\psi}_i^{k_i},$$
where $\overline{\psi}_i$ are the pullback of the first  Chern classes  of the tautological line bundles  over $\overline{M}_{0,n}(X,\beta)$ (which by abuse of notation we also denote by $\bar{\psi}_i$). 

For  classes $\delta_i\in H^*(\clG_{g_i},\mathbb{Q})$, 
set $\overline{\delta}_i=(\epsilon_{g_{i}}^{*})^{-1}(\delta_i)$. Descendant 
Gromov-Witten invariants  $\langle\overline{ \delta}_1\bar{\psi}_1^{k_1},\cdots,\overline{\delta}_n\bar{\psi}_n^{k_n}\rangle_{0, n, \beta}^X$ of $X$ are similarly defined. 
Theorem \ref{pushforward_vir} implies the following comparison result.
\begin{thm}\label{GW-inv1}
\begin{equation*}
 \langle \delta_1\bar{\psi}_1^{k_1},..., \delta_n\bar{\psi}_n^{k_n}\rangle_{0,n,\beta}^\clG=\frac{1}{r} \langle \overline{\delta}_1\bar{\psi}_1^{k_1},\cdots,\overline{\delta}_n\bar{\psi}_n^{k_n}\rangle_{0,n,\beta}^X.
\end{equation*}
Moreover, if $\vec{g}$ is not $\beta$-admissible, then the Gromov-Witten invariants of $\clG$ vanish.
\end{thm}

\begin{pf}
Denote by $\overline{ev}_i: \overline{M}_{0,n}(X, \beta)\to X$ the $i$-th evaluation map. Using the definition of $ev_i^*$ one can check that $ev_i^*(\delta_i)=p^*\overline{ev}_i^*(\overline{\delta}_i)$. Note also that $p^*\bar{\psi}_i=\bar{\psi}_i$. Thus using Theorem \ref{pushforward_vir} we have 
\begin{equation*}
\begin{split}
\langle \delta_1\bar{\psi}_1^{k_1},..., \delta_n\bar{\psi}_n^{k_n}\rangle_{0,n,\beta}^\clG&=\int_{[\mathcal{K}_{0,n}(\clG,\beta)^{\vec{g}}]^{vir}}\prod_{i=1}^n ev_i^*(\delta_i)\bar{\psi}_i^{k_i}\\ 
&=\int_{[\mathcal{K}_{0,n}(\clG,\beta)^{\vec{g}}]^{vir}}\prod_{i=1}^n p^*\overline{ev}_i^*(\overline{\delta}_i)\bar{\psi}_i^{k_i}\\
&=\int_{[\mathcal{K}_{0,n}(\clG,\beta)^{\vec{g}}]^{vir}}\prod_{i=1}^n p^*(\overline{ev}_i^*(\overline{\delta}_i)\bar{\psi}_i^{k_i})\\
&=\frac{1}{r}\int_{[\overline{M}_{0,n}(X,\beta)]^{vir}}\prod_{i=1}^n \overline{ev}_i^*(\overline{\delta}_i)\bar{\psi}_i^{k_i}\\
 &=\frac{1}{r} \langle \overline{\delta}_1\bar{\psi}_1^{k_1},\cdots,\overline{\delta}_n\bar{\psi}_n^{k_n}\rangle_{0,n,\beta}^X.
\end{split}
\end{equation*}
\end{pf}

In the following we use complex numbers $\cc$ as coefficients for the cohomology. For $\overline{\alpha}\in H^*(X, \mathbb{C})$ and an irreducible representation $\rho$ of $\mu_r$, we define $$\overline{\alpha}_\rho:=\frac{1}{r}\sum_{g\in \mu_r} \chi_{\rho}(g^{-1})\epsilon_g^*(\overline{\alpha}),$$ where $\chi_\rho$ is the character of $\rho$. The map $(\overline{\alpha},\rho)\mapsto \overline{\alpha}_\rho$  clearly defines an additive isomorphism
\begin{equation}\label{char_table_isom}
\bigoplus_{[\rho]\in \widehat{\mu_r}}H^*(X)_{[\rho]}\simeq H^*(I\clG, \mathbb{C}),
\end{equation}
where $\widehat{\mu_r}$ is the set of isomorphism classes of irreducible representations of $\mu_r$, and for $[\rho]\in\widehat{\mu_r}$ we define $H^*(X)_{[\rho]}:= H^*(X, \mathbb{C})$.

Theorem \ref{GW-inv1} together with orthogonality relations of characters of  $\mu_r$ implies the following 

\begin{thm}\label{GW-inv2}
Given $\overline{\alpha}_1, ...,\overline{\alpha}_n\in H^*(X, \bbQ)$ and integers $k_1,..., k_n\geq 0$, we have 
\begin{equation*}
\begin{split}
&\langle \overline{\alpha}_{1\rho_1}\bar{\psi}_1^{k_1},..., \overline{\alpha}_{n\rho_n}\bar{\psi}_n^{k_n} \rangle_{0,n,\beta}^{\clG}\\
=&\begin{cases}
\frac{1}{r^2}\langle \overline{\alpha}_1\bar{\psi}_1^{k_1},\cdots,\overline{\alpha}_n\bar{\psi}_n^{k_n}\rangle_{0,n,\beta}^X\chi_\rho(\exp(\frac{-2\pi \sqrt{-1} \int_\beta c_1(\mathcal{L})}{r})) &\text{if }\rho_1=\rho_2=...=\rho_n=:\rho,\\
0&\text{otherwise}\,.
\end{cases}
\end{split}
\end{equation*}
\end{thm}

\begin{pf}
By our definition we have 
\begin{equation*}
\begin{split}
&\langle \overline{\alpha}_{1\rho_1}\bar{\psi}_1^{k_1},..., \overline{\alpha}_{n\rho_n}\bar{\psi}_n^{k_n} \rangle_{0,n,\beta}^{\clG}\\
=&\frac{1}{r^n}\sum_{g_1,...,g_n\in \mu_r}\prod_{i=1}^n\chi_{\rho_i}(g_i^{-1})\langle \prod_{i=1}^n \epsilon_{g_i}^*(\overline{\alpha}_i)\bar{\psi}_i^{k_i} \rangle_{0,n,\beta}^\clG.
\end{split}
\end{equation*}
The term associated to $\vec{g}:=(g_1,...,g_n)$ in the above sum vanishes unless $\vec{g}$ is a $\beta$-admissible vector. This implies that $\prod_{i=1}^n g_i=\exp(\frac{2\pi \sqrt{-1}\int_\beta c_1(\mathcal{L})}{r})$. We rewrite this as $g_n^{-1}=\exp(\frac{-2\pi \sqrt{-1}\int_\beta c_1(\mathcal{L})}{r})\prod_{i=1}^{n-1} g_i$. Substitute this into above equation and use Theorem \ref{GW-inv1} to get 
\begin{equation*}
\begin{split}
&\langle \overline{\alpha}_{1\rho_1}\bar{\psi}_1^{k_1},..., \overline{\alpha}_{n\rho_n}\bar{\psi}_n^{k_n} \rangle_{0,n,\beta}^{\clG}\\
=&\frac{1}{r^n}\sum_{g_1,...,g_{n-1} \in \mu_r}\chi_{\rho_n}(\exp(\frac{-2\pi \sqrt{-1}\int_\beta c_1(\mathcal{L})}{r}))\left(\prod_{i=1}^{n-1}\chi_{\rho_i}(g_i^{-1})\chi_{\rho_n}(g_i)\right)\frac{1}{r}\langle \prod_{i=1}^n \overline{\alpha}_i\bar{\psi}_i^{k_i} \rangle_{0,n,\beta}^X.
\end{split}
\end{equation*}
Applying the orthogonality condition 
$$\frac{1}{r}\sum_{g\in \mu_r}\chi_\rho(g^{-1})\chi_{\rho'}(g)=\delta_{\rho, \rho'},$$
we find 
\begin{equation*}
\begin{split}
&\langle \overline{\alpha}_{1\rho_1}\bar{\psi}_1^{k_1},..., \overline{\alpha}_{n\rho_n}\bar{\psi}_n^{k_n} \rangle_{0,n,\beta}^{\clG}\\
=&\frac{1}{r}\chi_{\rho_n}(\exp(\frac{-2\pi \sqrt{-1}\int_\beta c_1(\mathcal{L})}{r}))\left(\prod_{i=1}^{n-1}\delta_{\rho_i, \rho_n}\right) \frac{1}{r}\langle \prod_{i=1}^n \overline{\alpha}_i\bar{\psi}_i^{k_i} \rangle_{0,n,\beta}^X.
\end{split}
\end{equation*}
The result follows.
\end{pf}

We now reformulate this in terms of generating functions. Let $$\{\overline{\phi}_i\, |\, 1\leq i\leq \text{rank}H^*(X, \mathbb{C})\}\subset H^*(X, \mathbb{C})$$ be an additive basis. According to the discussion above, the set $$\{\overline{\phi}_{i\rho}\,|\, 1\leq i\leq \text{rank}H^*(X, \mathbb{C}), [\rho]\in \widehat{\mu_r}\}$$ is an additive basis of $H^*(I\clG, \mathbb{C})$. Recall that the genus $0$ descendant potential of $\clG$ is defined to be 
\begin{eqnarray}
&\mathcal{F}^0_{\clG}(\{t_{i\rho, j}\}_{1\leq i\leq \text{rank}H^{*}(X,\mathbb{C}), \rho\in \widehat{\mu_r}, j\geq 0}; Q):= & \nonumber \\
& \sum_{\overset{n\geq 0, \beta\in H_2(X,\mathbb{Z})}{i_1,...,i_n; \rho_1,...,\rho_n; j_1,...,j_n}}\frac{Q^\beta}{n!}\prod_{k=1}^nt_{i_k\rho_k, j_k}\langle\prod_{k=1}^n \overline{\phi}_{i_k\rho_k}\bar{\psi}_k^{j_k} \rangle_{0,n,\beta}^{\clG}.& 
\end{eqnarray}
The descendant potential $\mathcal{F}^0_{\clG}$ is a formal power series in variables $t_{i\rho, j}, 1\leq i\leq \text{rank}H^*(X,\mathbb{C}), \rho\in \widehat{\mu_r}, j\geq 0$ with coefficients in the Novikov ring $\mathbb{C}[[\overline{NE}(X)]]$, where $\overline{NE}(X)$ is the effective Mori cone of the coarse moduli space of $\clG$. Here $Q^\beta$ are formal variables labeled by classes $\beta\in \overline{NE}(X)$. See e.g. \cite{ts} for more discussion on descendant potentials for orbifold Gromov-Witten theory. 

Similarly the genus $0$ descendant potential of $X$ is defined to be 
\begin{equation}\label{genus0_potential_of_X}
\mathcal{F}^0_X(\{t_{i,j}\}_{1\leq i\leq \text{rank}H^*(X,\mathbb{C}), j\geq 0};Q):=\sum_{\overset{n\geq 0, \beta\in H_2(X, \mathbb{Z})}{i_1,...,i_n; j_1,...,j_n}}\frac{Q^\beta}{n!}\prod_{k=1}^nt_{i_k, j_k}\langle \prod_{k=1}^n \overline{\phi}_{i_k}\bar{\psi}_k^{j_k}\rangle_{0,n,\beta}^X.
\end{equation}
The descendant potential $\mathcal{F}^0_X$ is a formal power series in variables $t_{i,j}, 1\leq i\leq \text{rank}H^*(X, \mathbb{C}), j\geq 0$ with coefficients in $\mathbb{C}[[\overline{NE}(X)]]$  and $Q^\beta$ is (again) a formal variable.
Theorem \ref{GW-inv2} may be restated as follows.

\begin{thm}\label{decomp_genus_0}
$$\mathcal{F}^0_{\clG}(\{t_{i\rho, j}\}_{1\leq i\leq \text{rank}H^*(X,\mathbb{C}), \rho\in \widehat{\mu_r}, j\geq 0}; Q)=\frac{1}{r^2}\sum_{[\rho]\in \widehat{\mu_r}}\mathcal{F}^0_X(\{t_{i\rho,j}\}_{1\leq i\leq \text{rank}H^*(X,\mathbb{C}), j\geq 0};Q_\rho),$$
where $Q_\rho$ is defined by the following rule: $$Q_\rho^\beta:=Q^\beta \chi_\rho\left(\exp\left(\frac{-2\pi \sqrt{-1} \int_\beta c_1(\mathcal{L})}{r}\right)\right),$$
and $\chi_\rho$ is the character associated to the representation $\rho$.
\end{thm}

Theorem \ref{decomp_genus_0} confirms the decomposition conjecture for genus $0$ Gromov-Witten theory of $\clG$.

We have another reformulation of Theorem \ref{decomp_genus_0}. Consider a new set of variables $$\{q_{i\rho, j}| 1\leq i\leq \text{rank}H^*(X, \mathbb{C}), [\rho]\in \widehat{\mu_r}, j\geq 0\}$$ 
defined by dilaton shifts
$$q_{i\rho,j}=\begin{cases}
t_{i\rho,j} & \text{ if } (i,j)\neq (1,1)\\
t_{1\rho,1}-1 & \text{ if } (i,j)=(1,1).
\end{cases}
$$
We view $\mathcal{F}_\clG^0$  as functions in the variables $q_{i\rho, j}$:

$$\mathcal{F}_\clG^0=\mathcal{F}_\clG^0(\{q_{i\rho,j}\}; Q).$$

Each term in the right-hand side of Theorem \ref{decomp_genus_0} can also be viewed as a function of the new variables $q_{i\rho,j}$ for a fixed $\rho\in \widehat{\mu_r}$:

$$\mathcal{F}_X^0(\{q_{i\rho, j}\}; Q_\rho).$$

By dilaton equation, we have $$\frac{1}{r^2}\mathcal{F}_X^0(\{q_{i\rho,j}\}; Q_\rho)=\mathcal{F}_X^0(\{\bar{q}_{i\rho, j}\}; Q_\rho)$$ where $\bar{q}_{i\rho, j}:= rq_{i\rho, j}$. Let $M_\rho$ denote the Frobenius structure on $H^*(X)_{[\rho]}=H^*(X)$ obtained using the potential function $\mathcal{F}_X^0(\{\bar{q}_{i\rho,j}\}; Q_\rho)$ and the Poincar\'e pairing on $X$. The following is a re-statement of Theorem \ref{decomp_genus_0}.

\begin{thm}\label{Frob_isom}
Under the isomorphism (\ref{char_table_isom}), the Frobenius structure defined by the genus $0$ Gromov-Witten theory of $\clG$ is isomorphic to $\oplus_{\rho} M_\rho$.
\end{thm}

\begin{rmk}
It is natural to ask for a generalization of Theorem \ref{decomp_genus_0} to higher genus Gromov-Witten theory. Suppose that the Frobenius structure associated to the genus $0$ Gromov-Witten theory of $X$ is generically semi-simple, then one can prove certain generalization of Theorem \ref{decomp_genus_0} to higher genus {\em ancestor invariants} by using Givental's formula \cite{gi}, \cite{teleman} to reduce the question to genus $0$. In \cite{AJT_root2} we will study the higher genus generalization of Theorem \ref{decomp_genus_0} in general (namely without assuming semi-simplicity). 
\end{rmk}

\appendix

\section{Banded abelian gerbes}\label{banded_abelian_gerbes_case}

Let $X$ be a smooth projective variety over $\bbC$. Let $G$ be a finite abelian group. The purpose of this Appendix is to explain (see Section \ref{ess-triv-Ggerbes-appendix}) how the results in the main part of the paper can be extended to banded $G$-gerbes $\clG$ over $X$ which are {\em essentially trivial}. 

We begin with some preliminary materials. Recall the following well-known  structure result for finite abelian groups: 
\begin{lem}
Let $G$ be a finite abelian group of order $N$. Then there exists a decomposition 
\begin{eqnarray}\label{G-decomp}
G\simeq \prod_{j=1}^k \mu_{r^{(j)}},\quad \text{where } \prod_{j=1}^k r^{(j)}=N. 
\end{eqnarray}
\end{lem}
Throughout this Appendix we fix such a decomposition (\ref{G-decomp}) of $G$. 

Observe that the inertia stack $I\clG$ admits a decomposition 
\begin{equation}\label{decomposition_of_inertia_stack}
I\clG=\cup_{g\in G} \clG_g,
\end{equation}
indexed by elements in $G$. Let $\bar{I}(\clG)_g\subset \bar{I}(\clG)$ be the image of $\clG_g$ under the natural map $I\clG\to \bar{I}(\clG)$ to the rigidified inertia stack. A vector of elements $$\vec{g}:=(g_1,...,g_n)\in G^{\times n}$$ is called {\em $\beta$-admissible} if the locus $$\clK_{0,n}(\clG, \beta)^{\vec{g}}:=\cap_{i=1}^n ev_i^{-1}(\bar{I}(\clG)_{g_i})$$
is non-empty. Note that for $1\leq i\leq n$ we may write $$g_i:=(g_i^{(1)},..., g_i^{(k)})\in \prod_{j=1}^k\mu_{r^{(j)}}=G.$$

\subsection{Essentially trivial abelian gerbes}\label{ess-triv-Ggerbes-appendix}

 By definition a $G$-gerbe over $X$ is {\em essentially trivial} if it  becomes trivial after contracted product with the trivial $\mathcal{O}_X^*$-gerbe. In this Section, let $\clG\to X$ is an essentially trivial $G$-banded gerbe over $X$. The following result is known (see e.g. \cite{FMN}, Proposition 6.9).
\begin{lem}\label{essentially_trivial=product_of_roots}
Let $\clG\to X$ is an essentially trivial $G$-banded gerbe over $X$, with $G$ finite and  abelian. Then there exist line bundles $\clL^{(1)}, ..., 
\clL^{(k)}$ over $X$ and positive integers $r^{(1)},..., r^{(k)}$, such that 
\begin{equation}\label{product_of_roots}
\clG\simeq \sqrt[r^{(1)}]{\clL^{(1)}/X}\times_X \sqrt[r^{(2)}]{\clL^{(2)}/X}\times_X...\times_X \sqrt[r^{(k)}]{\clL^{(k)}/X}.\end{equation}
\end{lem}

\begin{pf}
Let $[\clG]\in H^2_{et}(X, G)$ be the class of the gerbe $\clG$. Fix a decomposition  of $G$ as in (\ref{G-decomp}).  Denote by $p_j: G\to \mu_{r^{(j)}}$ the    projection to the $j$-th factor.  For $1\leq j\leq k$, the induced morphism 
$ p_{j*}:  H^2_{et}(X, G)\to   H^2_{et}(X, \mu_{r^{(j)}})$ maps the class of a $G$-banded 
gerbe $\clG$  to the class of the $\mu_{r^{(j)}}$-gerbe  
obtained  from $\clG$ by taking the contracted product with 
the trivial   $\mu_{r^{(j)}}$-gerbe . This is the same as taking the  rigidification of $\clG$ by  the  subgroup of the inertia $\overline{G}:=G/\mu_{r^{(j)}}$. The composition of $p_j$ with the standard embedding $\mu_{r^{(j)}}\to \bbC^*$ yields a homomorphism $\phi_j: G\to \bbC^*$. Clearly the composition $$G\overset{\phi_j}{\longrightarrow} \bbC^*\overset{(\cdot)^{r^{(j)}}}{\longrightarrow} \bbC^*$$ is trivial. 

Associated to the Kummer sequence $$1\longrightarrow \mu_{r^{(j)}}\longrightarrow \bbC^* \overset{(\cdot)^{r^{(j)}}}{\longrightarrow} \bbC^*\longrightarrow 1$$
there is a long exact sequence 
$$...\to \check{H}^1_{\acute{e}t}(X, \bbC^*)\to \check{H}^2_{\acute{e}t}(X, \mu_{r^{(j)}})\to \check{H}^2_{\acute{e}t}(X, \bbC^*)\to \check{H}^2_{\acute{e}t}(X, \bbC^*)\to ...$$
The map $\phi_j$ induces a homomorphism $\phi_{j*}: \check{H}^2_{\acute{e}t}(X, G)\to \check{H}^2_{\acute{e}t}(X, \bbC^*)$  mapping the class   $p_{j*}[\clG]$ of the   $\mu_{r^{(j)}}$-gerbe   obtained from $\clG$  by the homomorphism $p_j: G\to \mu_{r^{(j)}}$   to the class  of its contracted product with the trivial $\clO^*_X$-gerbe.  Since $\clG$ is essentially trivial, the class $\phi_{j*}([\clG])\in \check{H}^2_{\acute{e}t}(X, \bbC^*)$ is zero by definition. By the exact sequence above this means that there exists a line bundle $\clL^{(j)}$ over $X$ such that the $\mu_{r^{(j)}}$-gerbe of class  $p_{j*}[\clG]$  is isomorphic to the root gerbe $\sqrt[r^{(j)}]{\clL^{(j)}/X}$. 
We can prove the claim by induction on the number $k$  of cyclic groups appearing in the decomposition of $G$.   For $k=1$ the claim is true by definition. 
Assume it is true for $k=n-1$. Let $\clG$ be a $G$-banded gerbe, where 
$G\simeq \prod_{j=1}^n \mu_{r^{(j)}}$.  Consider the group homomophisms 
 induced by the  rigidification 
$$ \check{H}^2_{\acute{e}t}(X, G)\to   \check{H}^2_{\acute{e}t}(X, \overline{G}) \oplus \check{H}^2_{\acute{e}t}(X, \mu_{r^{(n)}}),   $$
where $\overline{G}\simeq G/\mu_{r^{(n)}}$. We denote the corresponding gerbes by
$\overline{\clG}$ and $\overline{\clG}_k$. We have a commutative diagram
\begin{eqnarray}
\xymatrix{
\clG\ar@{-->}[rd]\ar@/^1pc/[rrd]\ar@/_1pc/[ddr] & & \\
&  \tilde{\clG} \ar[d]\ar[r]\ar@{}[rd]|{\square} & \overline{\clG}\ar[d]\\
& \overline{\clG}_k\ar[r] & X, 
}
\end{eqnarray}  
where the dotted arrow is induced by the universal property of the fiber product. The morphism  $\clG\to \overline{\clG} \times \overline{\clG}_k$ is representable and factors through $\clG\to \tilde{\clG}$, which is therefore representable. We conclude by observing that a representable morphism between two gerbes banded by the same group is an isomorphism. 
\end{pf}

In view of Lemma \ref{essentially_trivial=product_of_roots} we assume that the gerbe $\clG$ is of the form (\ref{product_of_roots}). We call this gerbe a {\em multi-root gerbe}. Recall that to give a morphism $Y\to \clG$ is the same as giving a morphism $f: Y\to X$ and line bundles $M_1, ..., M_k$ over $Y$ together with isomorphisms $\phi_j: M_j^{\otimes r^{(j)}}\simeq f^*\clL^{(j)}$, $1\leq j\leq k$. 

The constructions in Section \ref{moduli_maps_root_gerbe} can be easily modified to treat the multi-root gerbe $\clG$.  Arguments proving Lemma \ref{adm_vect_lemm} and Proposition \ref{unique_geom_lift} easily yield the following 
\begin{prop}\label{adm-vect-geom-lift-Ggerbe-prop}
\hfill
\begin{enumerate}
\item
A vector $\vec{g}$ is $\beta$-admissible (with respect to a class $\beta\in H_2^+(X, \bbZ)$) if and only if 
\begin{equation}\label{phase_shift_eq_multi_roots}
\prod_{i=1}^n g_i^{(j)}=\exp\left(\frac{2\pi \sqrt{-1}}{r^{(j)}}\int_\beta c_1(\clL^{(j)})\right),\quad  1\leq j\leq k.
\end{equation}
\item
Given a vector $\vec{g}$ satisfying (\ref{phase_shift_eq_multi_roots}) and a stable map $[f:(C, p_1,...,p_n) \to X]\in \overline{M}_{0,n}(X,\beta)(\bbC)$, there exists, up to isomorphisms, a unique twisted stable map $\tilde{f}:(\clC, \sigma_1,..., \sigma_n)\to \clG$ in $\clK_{0,n}(\clG,\beta)^{\vec{g}}$ lifting $f$.
\end{enumerate}
\end{prop}

\begin{rmk}
The $n$-tuple $(g_1^{(j)},...,g_n^{(j)})$ is a $\beta$-admissible vector for the root gerbe $\sqrt[r^{(j)}]{\clL^{(j)}/X}$ as in Definition \ref{adm_vector-def}.
\end{rmk}

Next we define some numbers. 
\begin{defn}
\begin{enumerate}
\item
For $1\leq i\leq n$, let $r_i$ be the order of $g_i$ in $G$. Each $g^{(j)}_i, 1\leq i\leq n$ may be identified with a root of unity $$g^{(j)}_i=\exp(2\pi \sqrt{-1} \theta^{(j)}_i), \quad  \text{where } \theta^{(j)}_i\in \mathbb{Q}\cap [0, 1),$$
which defines the rational numbers $\theta^{(j)}_i, 1\leq i\leq n$. For $1\leq i\leq n$ and $1\leq j\leq k$, define 
\begin{equation}\label{the_triple_of_integers_multi_root}
\rho^{(j)}_i:= r^{(j)}\theta^{(j)}_i, \quad r^{(j)}_i:=\frac{r^{(j)}}{gcd (r^{(j)}, \rho^{(j)}_i)}, \quad m^{(j)}_i:=\frac{\rho^{(j)}_i}{gcd (r^{(j)}, \rho^{(j)}_i)}.
\end{equation}
 Note that $r_i^{(j)}$ divides $r_i$, and $r_i^{(j)}$ is the order of $g_i^{(j)}$ in $\mu_{r^{(j)}}$.
\item
For a pair $(T, \beta')$ indexing the boundary divisors of $\MM_{0,n,\beta}$ as in Definition \ref{boundary_of_M0nbeta}, define 
\begin{equation}\label{triple_of_integers_for_node_multi_root}
\theta^{(j)}_{T, \beta'}:=\langle\frac{1}{r^{(j)}}\int_{\beta'} c_1(\mathcal{L}^{(j)})-\sum_{i\in T} \theta^{(j)}_i\rangle, \quad r^{(j)}_{T, \beta'}:=\frac{r^{(j)}}{gcd(r^{(j)}, r^{(j)}\theta^{(j)}_{T, \beta'})}, \quad m^{(j)}_{T, \beta'}:=\frac{r^{(j)}\theta^{(j)}_{T, \beta'}}{gcd(r^{(j)}, r^{(j)}\theta^{(j)}_{T, \beta'})}.  
\end{equation}
\item
Define $$g_{T, \beta'}^{(j)}:=\exp(2\pi\sqrt{-1}\theta_{T, \beta'}^{(j)})\in \mu_{r^{(j)}}, \quad g_{T, \beta'}:=(g_{T, \beta'}^{(1)},..., g_{T,\beta'}^{(k)})\in G.$$ And let $r_{T, \beta'}$ be the order of $g_{T, \beta'}$ in $G$. 
\end{enumerate}
\end{defn}

With the numbers defined above, the constructions and results in Sections \ref{section:stack_YY} and \ref{section:gerbe_structure} are valid for the multi-root gerbe $\clG$. The proofs are straightforward modifications. In particular we still have the diagram (\ref{diag_root_gerb_mod_sp}). 

Moreover Theorem \ref{mod_space_to_root_gerbe_is_root_gerbe} admits a generalization to multi-root gerbes:

\begin{thm}\label{mod_space_to_multi-root_gerbe_is_multi-root_gerbe}
$\clK_{0,n}(\clG, \beta)^{\vec{g}}$ is a multi-root gerbe over $P_n^{\vec{g}}$. 
\end{thm}

To prove Theorem \ref{mod_space_to_multi-root_gerbe_is_multi-root_gerbe} it suffices to repeat the arguments in the proof of Theorem \ref{mod_space_to_root_gerbe_is_root_gerbe} multiple times. The key point is to construct a collection of line bundles, generalizing the one in (\ref{L_YY_eq}):
\begin{equation}\label{L_YY_eq_j}
\clL_\YY^{(j)}:=\mathcal{O}_{\YY_{0,n+1, \beta}^{\vec{g}\cup \{1\}}}\left(\sum_{1\leq i\leq n} \frac{d^{(j)}_i}{r^{(j)}_i} S_i -\sum_{(T, \beta')\in \mathcal{I}_D} \frac{d^{(j)}_{T, \beta'}}{r^{(j)}_{T, \beta'}}D_{\beta'}^{T\cup\{n+1\}}\right), \quad 1\leq j\leq k.
\end{equation}
Here the set $\mathcal{I}_D$ is defined on page \pageref{set_I_sub_D}, and we use the following definition:
\begin{defn}
\begin{enumerate}
\item
For $1\leq i\leq n$ and $1\leq j\leq k$ we define $d^{(j)}_i\in \bbZ$ by requiring 
\begin{equation}\label{requirement_for_d_i_j}
g^{(j)}_i=\exp(\frac{2\pi \sqrt{-1}}{r^{(j)}_i}d^{(j)}_i), \quad \text{and } \sum_{i=1}^n \frac{d^{(j)}_i}{r^{(j)}_i}=\frac{1}{r^{(j)}}\int_{\beta}c_1(\clL^{(j)}).
\end{equation}

\item
To a pair $(T, \beta')$ which indexes a boundary divisor of $\MM_{0,n,\beta}$, we associate integers $d^{(j)}_{T, \beta'}$ such that 
\begin{equation}\label{requirement_for_d_node_j}
\sum_{i\in T} \frac{d^{(j)}_i}{r^{(j)}_i} +\frac{d^{(j)}_{T, \beta'}}{r^{(j)}_{T, \beta'}}=\frac{1}{r^{(j)}}\int_{\beta'} c_1(\clL^{(j)}), \quad 1\leq j\leq k.
\end{equation}
\end{enumerate}
\end{defn}


 \end{document}